\documentclass[12pt]{article}

\usepackage{amsmath}
\usepackage{accents}
\usepackage{amssymb}
\usepackage{epsfig}
\usepackage{xspace}
\usepackage{float} 
\usepackage{multirow}
\usepackage{pgf}
\usepackage{soul}
\usepackage[top=2cm, bottom=2cm, left=2cm, right=2cm]{geometry}
\usepackage{minitoc}
\usepackage{bbold}
\usepackage[nottoc,numbib]{tocbibind}
\usepackage{xstring}
\usepackage{hyperref}

\usepackage{enumitem}
\setlist{noitemsep}
\usepackage{algorithm}

\definecolor{Red}{rgb}{1,0,0}
\definecolor{Green}{rgb}{0,0.6,0.4}
\definecolor{Blue}{rgb}{0,0,1}
\definecolor{Gray}{rgb}{0.2,0.2,0.2}

\def\R{{\mathbb{R}}}
\def\N{{\mathbb{N}}}
\def\P{{\mathbb{P}}}
\def\E{{\mathbb{E}}}

\newcommand{\y}{\mathbf{y}}

\newcommand{\x}{\mathbf{x}}
\newcommand{\z}{\mathbf{z}}
\newcommand{\Z}{\mathbf{Z}}

\newcommand{\W}{\mathbf{W}}

\renewcommand{\u}{\mathbf{u}}

\newcommand{\mads}{{MADS}\xspace}
\renewcommand{\emptyset}{\varnothing}

\newcommand{\fh}{\hat{f}}
\newcommand{\ch}{\hat{c}}

\newcommand{\yh}{\hat{y}}
\newcommand{\hh}{\hat{h}}

\newlength{\dhatheight}
\newcommand{\doublehat}[1]{
    \settoheight{\dhatheight}{\ensuremath{\hat{#1}}}
    \addtolength{\dhatheight}{-0.35ex}
    \hat{\vphantom{\rule{1pt}{\dhatheight}}
    \smash{\hat{#1}}}
}

\newlength{\dwidehatheight}
\newcommand{\widedoublehat}[1]{
    \settoheight{\dwidehatheight}{\ensuremath{\hat{#1}}}
    \addtolength{\dwidehatheight}{-0.35ex}
    \widehat{\vphantom{\rule{1pt}{\dwidehatheight}}
    \smash{\widehat{#1}}}
}

\newcommand{\cvvalue}[1]{\doublehat{#1}}
\newcommand{\ycv}{\cvvalue{y}}
\newcommand{\hcv}{\cvvalue{h}}
\newcommand{\ccv}{\cvvalue{c}}
\newcommand{\fcv}{\cvvalue{f}}

\newcommand{\prech}{\;\widehat{\prec}\;}
\newcommand{\preceqh}{\;\widehat{\preceq}\;}
\newcommand{\preccv}{~\widedoublehat{\prec}~}

\newcommand{\xib}{{\boldsymbol{\xi}}}
\renewcommand{\t}{^\top}

\newcommand{\sinf}{\s^\infty}
\newcommand{\sbest}{\s^*}
\newcommand{\XUS}{\X\cup\S}
\newcommand{\AND}{\xspace~{\tt and}~\xspace}
\newcommand{\xor}{\xspace~{\tt xor}~\xspace}

\renewcommand{\chi}{\mathcal{X}}

\newcommand{\X}{\mathbf{X}}
\renewcommand{\P}{\mathbf{P}}
\newcommand{\Xh}{\mathbf{\hat{X}}}

\newcommand{\ind}{\mathbb{1}}

\newcommand{\metric}{\mathcal{E}}
\mathchardef\mhyphen="2D

\newcommand{\mesh}{\mathcal{M}}

\newcommand{\card}{\text{card}}

\newcounter{modelIndex}
\addtocounter{modelIndex}{1}

\renewcommand{\quote}[1]{``#1''}
\newcommand{\figWidth}{12cm}

\newcommand{\ignore}[1]{}

\renewcommand{\t}{^\top}

\newcommand{\Y}{\mathbf{Y}}
\renewcommand{\S}{\mathbf{S}}

\newcommand{\e}{\mathbf{e}}

\newcommand{\s}{\mathbf{s}}
\renewcommand{\u}{\mathbf{u}}

\newcommand{\xlb}{\underline{\x}}
\newcommand{\xub}{\bar{\x}}
\newcommand{\xstar}{\x^*}

\title{Parallel Surrogate-assisted Optimization Using Mesh Adaptive Direct Search}
\author{Bastien Talgorn$^1$, St\'ephane Alarie$^2$, and Michael Kokkolaras$^1$}

\date{
    \small{
    $^1$McGill University, GERAD, Montr\'eal, Qu\'ebec, Canada\\%
    $^2$Hydro-Qu\'ebec's Research Institute, GERAD, Montr\'eal, Qu\'ebec, Canada\\%
    }
}

\begin{document}

\maketitle

\abstract{
We consider computationally expensive blackbox optimization problems and present a method that employs surrogate models and concurrent computing at the search step of the mesh adaptive direct search (MADS) algorithm. 
Specifically, we solve a surrogate optimization problem using locally weighted scatterplot smoothing (LOWESS) models to find promising candidate points to be evaluated by the blackboxes. We consider several methods for selecting promising points from a large number of points. 
We conduct numerical experiments to assess the performance of the modified MADS algorithm with respect to available CPU resources by means of five engineering design problems.}

\setlength{\parindent}{0pt}

\section{Introduction}
\label{sec:intro}

We consider the optimization problem
	\begin{equation}
		\begin{array}{rl}
			\underset{\x \in \chi}{\min}	&   f(\x)\\
			\text{subject to}			&   c_j(\x) \le 0,  ~~ j = 1, 2, \hdots, m,
		\end{array}
		\label{eq:MainProblem}
		\tag{$P$}
	\end{equation}
	where $f(\x)$ is the objective function,
	$\x \in \R^n$ is the vector of decision variables, 
	$\chi$ is a subset of $\R^n$, 
	and $c_j(\x)$ are general nonlinear constraints.
We assume that some (at least one) of the functions  $\{f,  c_1, c_2, \hdots, c_m\}$
	are evaluated using simulations
	or other computational procedures that are blackboxes.
In particular,
	we consider the case where these blackboxes are computationally expensive,
	possibly nonsmooth and/or nonconvex, and that the process used to evaluate them
	may crash or fail to return a value.
Finally, we assume that function gradients either do not exist theoretically or, if they do, cannot be computed  
	or approximated with reasonable computational effort.

Metaheuristics and derivative-free search algorithms are commonly used for solving~\eqref{eq:MainProblem}.
The former
	(e.g., genetic algorithms(GAs), particle swarm optimization~(PSO), tabu search~(TS), etc.)
	are commonly used for global exploration
	while the latter
	(e.g., generalized pattern search~(GPS), mesh adaptive direct search~(MADS), and trust-region methods~(DFO, COBYLA, CONDOR))
	are local methods with convergence properties~\cite{DBLP:series/sci/KramerCK11}. 
In this work, we use the NOMAD implementation~\cite{AuCo04a, Le09b} of the MADS algorithm~\cite{AuDe2006} to solve~\eqref{eq:MainProblem}.

Multiprocessor computers,
	supercomputers,
	cloud computing,
	or just a few connected PCs
	can provide parallel (or concurrent) computing opportunities to speed up so-called trajectory-based optimization algorithms.
According to \cite{AlbaLuqueNesmachnow2013}, 
	three ways are commonly used to achieve this:
		(i)~parallel evaluation of neighborhood solutions (distributed evaluations),
		(ii)~parallel trajectories from the same (or different) initial guess(es) (independent optimization runs),
		(iii)~the evaluation of a point $\x$ is performed in parallel (i.e., the search is sequential).
The implementation of (iii) depends only on the blackbox, while the other two are related to the optimization algorithm.

NOMAD offers implementations for (i) and (ii) through p-MADS for (i) and Coop-MADS for (ii)~\cite{LeAbAuDe10}.
In both cases,
	the parallel evaluations are handled by NOMAD by means of MPI calls to the blackbox~\cite{Le09b}.
However,
	if one prefers to take charge of the distribution of evaluations,
	one can implement p-MADS with blocks by using NOMAD in batch mode~~\cite{LeAbAuDe10}.
In that case, instead of using MPI calls,
	NOMAD writes all points to be evaluated in an input file,
 waits for all evaluations to be completed, 
	and reads the obtained values of $f(\x)$ and $c_j(\x)$ from an output file.
Then, NOMAD 
either stops if a local optimum is reached
	or submits a new block of points to evaluate.
	We use here the block functionality of NOMAD,  
	adopting option (i) for parallel evaluations, because of its flexibility and generality.

The MADS algorithm includes two steps at each iteration, the SEARCH and the POLL.
The SEARCH step is flexible (defined by the user)
	and aims at determining one or more new points $\x \in \chi$ that improves the current best solution.
The POLL step is defined
	according to the convergence analysis of the algorithm and
generates trial points around the current best solution.
The number of poll points at each iteration is either $2n$ or $n+1$ depending on the utilized pattern with
$n$ being the size of $\x$.

The number of points evaluated in the SEARCH step depends on the methods chosen or defined by the user.
Several techniques are already implemented and available in NOMAD, including
	the speculative search (SS),
	the quadratic model search (QUAD),
	the variable neighborhood search (VNS),
	the Latin hypercube search (LHS), and
	the Nelder Mead search (NMS).
One can also implement their own technique if one so desires, 
	which is called a \emph{user search} (US).

With the exception of LHS,
all provided techniques usually return only one trial point.
When several techniques are used at once,
	they are called one after the other along the SEARCH step,
	each technique providing its own trial point,
	which is evaluated by the blackbox before proceeding to the next technique.
Assuming that $2n$ CPUs are available for solving \eqref{eq:MainProblem},
	the POLL step can make good use of these CPUs. 
	However, since SEARCH step
	evaluations are sequential,
 progress is slowed down with almost all CPUs being idle.
One may argue that we should then only use LHS since it can generate $2n$ points.
However, since LHS is random,
	its points will quickly become less promising after a few iterations.

Considering that the number of available CPUs are now, particularly with the emergence of cloud computing, relatively inexpensive and unlimited,
	we should rethink the SEARCH step to be as effective as the POLL step in terms of CPU use.
In this work,
	we propose a SEARCH step technique that returns a large number of diverse points for evaluation. 

The paper is structured as follows.
The general idea behind the proposed technique is described in Section~\ref{sec:genIdea}.
In Section~\ref{sec:parallel},
	six different methods are presented for selecting various candidates from a large set of points.
In Section~\ref{sec:algo},
	the resulting model search for parallel computing is specified.
In Section~\ref{sec:test}, 
	we test our SEARCH step technique on five engineering design optimization problems using up to 64 processors.
A discussion concludes the paper.

\section{Proposed SEACH step technique}
\label{sec:genIdea}

One of the practical challenges of the SEARCH step	is that only one candidate is obtained at a significant computational investment~\cite{TaLeDKo2014, PoTaHaKoLe2014, Ensemble2016, Lowess2016}.
Specifically, regardless of the number of available CPUs,
	only one CPU is used in the SEARCH step for blackbox evaluations, with the exception of LHS.
	Before presenting our idea for mitigating this practical challenge, we will assume that computationally inexpensive surrogate models of the expensive blackboxes are available.
We can then consider the surrogate problem of problem~\eqref{eq:MainProblem}
	\begin{equation}
		\begin{array}{rl}
			\underset{\x \in \chi}{\min}	&   \fh(\x)\\
			\text{subject to}			&   \ch_j(\x) \le 0,  ~~ j = 1, 2, \hdots, m,
		\end{array}
		\label{eq:SurrogateProblem}
		\tag{$\hat{P}$}
	\end{equation}
	where  $\{\fh, \ch_1, \ch_2, \hdots, \ch_m\}$ are surrogate models of $\{f,  c_1, c_2, \hdots, c_m\}$, respectively.
Note that we only need to ensure that the minimizers of~\eqref{eq:MainProblem} and~\eqref{eq:SurrogateProblem} are close enough, and not that the surrogate models are good approximations of the blackboxes globally. 
It then follows that
	a minimizer of~\eqref{eq:SurrogateProblem} will be a good candidate for the solution of~\eqref{eq:MainProblem}.

If both problems have the same minimizers, they may share features in other areas of~$\chi$ as well.
Since the evaluations of  $\fh(\x)$ and $\ch_j(\x)$ are rather inexpensive
	compared to $f$ and $c_j$,
	one can allow a very large budget of model evaluations
	to solve  \eqref{eq:SurrogateProblem},
	extending thus the number of design space areas that will be visited.
This is acceptable as long as the solution of \eqref{eq:SurrogateProblem} is faster
	than any single evaluation of the blackboxes. 
Considering there are $q$ CPUs available for blackbox evaluations,
	one may then select $q$ points from the available budget
	by solving \eqref{eq:SurrogateProblem}.
The $q$ points can be selected to consider areas of~$\chi$ 
	that have been neglected until now in the solution of~\eqref{eq:MainProblem}.

The above proposition proposes the use of surrogate models
	$\{\fh, \ch_1, \ch_2, \hdots, \ch_m\}$ in a manner 
	that is not reported in~\cite{Haftka2016}, which mentions two ways of exploiting surrogates in the context of parallelization.
The simplest is to fit $q$ different surrogate models
	at the same points already evaluated by the blackbox functions.
This allows to get $q$ different promising candidates
	and requires no uncertainty quantification for the surrogate models.
One can also combine the surrogates; 
distance-based criteria can be added to ensure diversity between the candidates.
The other way is to use a single surrogate model 
	and consider $q$ points where the blackboxes should be evaluated at
	to improve its accuracy.
%

We propose an intermediate approach.
We use only one surrogate model for each blackbox. 
The $q$ candidates are extracted from that single surrogate, but not with the aim of improving it.
Instead, the $q$ candidates are selected to be the most interesting to advance the optimization process.
%

\section{Methods for selecting candidate points}
\label{sec:parallel}

Let $\X = \{\x_1 , \x_2, \hdots, \x_k\}$ be the set of all points evaluated by the blackbox.
Note that $\X \subset \chi \subset \R^n$.
We denote $\Xh$ the \emph{surrogate cache},
	i.e., the set of all points for which $\{\fh, \ch_1, \ch_2, \hdots, \ch_m\}$ have been evaluated
	during the solution of~\eqref{eq:SurrogateProblem}.
Similarly, $\Xh \subset \chi \subset \R^n$.
Let $\S$ be the set of  points that are selected by the SEARCH step 
	to be evaluated with the blackbox.
The set $\S$ is initially empty 
	and is built from the points of $\Xh$
	(ensuring that $\S \subset \Xh $)
	with a greedy algorithm by means of up to six  selection methods,
each one having a different goal.
	\begin{itemize}
		\item Method~\ref{algo:method1} selects the best point of $\Xh$  not in $\XUS$;
		\item Method~\ref{algo:method2} selects the most distant point of $\Xh$ from $\XUS$;
		\item Method~\ref{algo:method3} selects the best point of $\Xh$ at a certain distance of $\XUS$;
		\item Method~\ref{algo:method4} selects the best point of $\Xh$ under additional constraints;
		\item Method~\ref{algo:method5} selects a point of $\Xh$ that is a possible local minimum of the surrogate problem; 
		\item Method~\ref{algo:method6} selects a point of $\Xh$ in a non-explored area.
	\end{itemize}
Note that some selection methods may fail to return a candidate,
	particularly methods~\ref{algo:method3}  and \ref{algo:method4}.
If this happens, the next method is used.
We repeat and loop through all methods until we obtain enough points in $\S$ matching the available CPUs.
Some points $\x \in \Xh$ can also belong to the set $\X$; this is not an issue since all methods only select point $\x \in \Xh$ to be added to $\S$
	if and only if $\x \notin \X$.
The selection methods are detailed below after some definitions.

\subsection{Definitions}

Let $d(A,B)$ be the Euclidean distance between two subsets $A$ and $B$ of $\R^n$
	\begin{equation}
		d(A,B) = \underset{a \in A}{\min} \; \underset{b \in B}{\min} \; \|a-b\|_2.
	\end{equation}
As a convention,
	the distance to an empty set is infinite:
	$d(A,\emptyset) = d(\emptyset,\emptyset) = +\infty$.
By extension,
	we will denote the distance between an element $a \notin B$ and the subset $B$
	simply by $d(a,B)$,
	which implies that $a$ also refers to the particular subset containing only $a$,
	i.e., $\{a\}$.

Regarding feasibility, we consider the aggregate constraint violation function used in~\cite{FlLe02a},
	i.e.,
	$ h(\x) = \sum_{j=1}^m \max\{0, c_j(\x)\}^2 $.
The same function is used in the \emph{progressive barrier} mechanism in NOMAD~\cite{AuDe09a}.

We also define the order operators between two points $\x$ and $\x' \in \chi$:
	\begin{align}
		\x \prec \x' \Leftrightarrow & \left\{
			\begin{array}{l}
				h(\x) < h(\x')\\
				\text{or}\\
				h(\x) = h (\x') \text{ and } f(\x) < f(\x') 
 			\end{array}
			\label{eq:order}
		\right.\\
		\x \preceq \x' \Leftrightarrow & \;\; {\tt not}(\x' \prec \x),
		\label{eq:ordereq}
	\end{align}
	which are transitive.
By those definitions,
	an incumbent solution $\x^\prime$ of the original problem~\eqref{eq:MainProblem}
	is such that $\x^\prime \preceq \x,\; \forall \x \in \X$.
Similarly,
	a global minimizer $\x^*$ is such that $\x^*\preceq \x,\; \forall \x \in \chi$.
In the same manner as we define $\prec$ and $\preceq$ for $f$ and $h$,
	we define $\prech$ and $\preceqh$ for $\fh$ and $\hh$. 
%

Finally,
	to simplify the description of the proposed selection methods,
	 we define $\sinf$ as a \emph{virtual point}
	(in the sense that it does not have coordinates in $\R^n$), 
	which represents the worst possible candidate in $\R^n$:
	\begin{equation}
		\hh(\sinf) = \fh(\sinf) = +\infty \text{ \;\; and \;\; } d(\sinf,\X)=0.
		\label{eq:sinf}
	\end{equation}

\subsection{Detailed description of selection methods}
\label{sec:detailedDescription}

\newenvironment{myAlgo} [1][]
{
\setlength\parindent{0pt}
\hspace*{0.1\linewidth}\begin{minipage}{0.8\linewidth}
\begin{algorithm}[H]}
{\end{algorithm}
\end{minipage}\\\\\\
}

\newcommand{\inArray}[1]{\begin{array}{l}#1\end{array}}
\newcommand{\algoNew}[1]{
	$\begin{array}{l}
		#1
	\end{array}$
}
\newcommand{\algoLine}[1]{\text{#1}\\}
\newcommand{\algoSection}[2]{\text{\textbf{#1}}\\
	\;\;\;\;\begin{array}{l}
			#2
	\end{array}\\
}
\newcommand{\algoOpen}[2]{
	\algoLine{#1}
	\left|\;
		\begin{array}{l}
			#2
		\end{array}
	\right.	\\
}
\newcommand{\algoFont}[1]{\textbf{\tt #1}}
\newcommand{\algoAfterBlockVSpace}{\vspace{0.08cm}} 
\newcommand{\algoBlockVerticalLine}{\hspace{0.3cm}} 
\newcommand{\algoBlockVerticalSpace}{\hspace{-0.1cm}} 

\newcommand{\algoBox}[3]{
	\begin{array}{| l @{} l}
			\cline{1-1}
			\multicolumn{2}{|l}{\algoBlockVerticalSpace\text{#1}  }\\
			\algoBlockVerticalLine & \inArray{#2} \\
			\multicolumn{2}{|l}{\algoBlockVerticalSpace\text{#3} }\\
			\cline{1-1}
	\end{array}\algoAfterBlockVSpace\\
}

\newcommand{\algoBoxIfElse}[5]{
	\begin{array}{| l @{} l}
			\cline{1-1}
			\multicolumn{2}{|l}{\algoBlockVerticalSpace\text{#1}  }\\
			\algoBlockVerticalLine & \inArray{#2} \\
			\multicolumn{2}{|l}{\algoBlockVerticalSpace\text{#3} }\\
			\algoBlockVerticalLine & \inArray{#4} \\
			\multicolumn{2}{|l}{\algoBlockVerticalSpace\text{#5}  }\\
			\cline{1-1}
	\end{array}\algoAfterBlockVSpace\\
}
\newcommand{\algoEnd}{\algoLine{end}}


{\bf Method~\ref{algo:method1}.}
The first selection method selects the best point $\s$ of $\Xh$ under the constraint that $d(\s,\XUS) > 0$,
	which means that $\s$ is not in the set $\X$ of evaluated points nor already selected
	(i.e., $\notin \S$).
This method reflects how surrogate models are typically used for finding new candidate points.

\begin{myAlgo}[!ht]
	\algoNew
	{
		\algoLine{$\sbest \leftarrow \sinf$}

		\algoOpen{for all $\s \in \Xh$, do:}
		{
			\algoOpen{if $\s \prech \sbest$ \AND $d(\s,\XUS) > 0$, then:}
			{
				\algoLine{$\sbest \leftarrow \s$}
			}
 			\algoEnd
		}
		\algoEnd

		\algoOpen{if $\sbest \neq \sinf$, then:}
		{
			\algoLine{$\S \leftarrow \S \cup \{\sbest\}$}
		}
		\algoEnd
	}
	\caption{Selection of the best point (Method~\ref{algo:method1})}
	\label{algo:method1}
\end{myAlgo}	


{\bf Method~\ref{algo:method2}.}
The second method aims to maximize the diversity of the candidates to be evaluated.
It selects the point $\s$ of $\Xh$ that maximizes the distance $d(\s,\XUS)$,
	i.e., as far as possible from points already evaluated.
	
%
%
\begin{myAlgo}[!ht]
	\algoNew
	{ 
 		\algoLine{$\sbest \leftarrow \sinf$}

		\algoOpen{for all $\s \in \Xh$, do:}
		{
			\algoOpen{if $d(\s,\XUS) > d(\sbest,\XUS)$, then:}
			{
				\algoLine{$\sbest \leftarrow \s$}
			}
			\algoEnd
		}
		\algoEnd

		\algoOpen{if $\sbest \neq \sinf$, then:}
		{
			\algoLine{$\S \leftarrow \S \cup \{\sbest\}$}
		}
		\algoEnd
	}
	\caption{Selection of the most distant point to $\XUS$ (Method~\ref{algo:method2})}
	\label{algo:method2}
\end{myAlgo}

 
{\bf Method~\ref{algo:method3}.}
This method selects the best point $\s$ of $\Xh$ under the constraint that $d(\s,\XUS) \ge d_{\min}$,
	where $d_{\min}$ is initialized at 0 at the beginning of the selection process
	and increased progressively as the method is applied. 
Method~\ref{algo:method3} may fail to select a candidate $\s$ when $d_{\min}$ becomes too large.
Since the selected points $\S$ must be projected on the current mesh
 $ \mesh = \{\x+\Delta^m \mathbf{D} \z, \z\in \N^{n_D}, \x \in \X\}$
as required by MADS,
	incrementing $d_{\min}$ by the current mesh size $\Delta^\mesh$
	allows to avoid that several candidates become identical after the projection.  

\begin{myAlgo}[!ht]
	\algoNew
	{
		\algoLine{$\sbest \leftarrow \sinf$}

		\algoOpen{if first use of Method~\ref{algo:method3}, then:}
		{
			\algoLine{$d_{\min} \leftarrow 0$}
		}
		\algoEnd  

		\algoOpen{for all $\s \in \Xh$, do:}
		{
			\algoOpen{if $\s \prech \sbest$ \AND $d(\s,\XUS) \ge d_{\min}$, then:}
			{
				\algoLine{$\sbest \leftarrow \s$}
			}
			\algoEnd  
		}
		\algoEnd

		\algoOpen{if $\sbest \neq \sinf$, then:}
		{
			\algoLine{$\S \leftarrow \S \cup \{\sbest\}$}
			\algoLine{$d_{\min} \leftarrow d_{\min} + \Delta^\mesh$}
		}
		\algoEnd
	}
	\caption{Selection of the best point with a constraint on the distance to $\XUS$ (Method~\ref{algo:method3})}
	\label{algo:method3}
\end{myAlgo}


\newcommand{\chmax}{\ch_{\max}}
\newcommand{\chmargin}{\ch_{\text{margin}}}

{\bf Method~\ref{algo:method4}.}
Considering that the surrogate models $\ch_j$ may fail to predict correctly if $\s$ is feasible,
	the present method tries to select points
	that will be likely to be feasible when evaluated by the blackboxes $c_j$.
This is done by selecting the best feasible point $\s$ of $\Xh$
	under the constraint  $\chmax(\s) \le \chmargin$, 
	where $\chmax(\s)$ is defined as being the most violated constraint of $\s$,
	i.e.,
	\begin{equation}
		\chmax(\s) = \underset{j=1,2,\dots,m}{\max}\; \ch_j(\s),
	\end{equation}
	where $\chmargin$ is set  as
	\begin{equation}
		\chmargin \leftarrow \underset{\substack{\s\in\Xh\\ \chmax(\s)<0}}{\max}\chmax(\s)
	\end{equation}
	and quantifies,
	among all feasible points of $\Xh$,
	the smallest amount by which these are satisfied.
	
By definition,
	$\chmax(\s) \le 0$ if $\s$ is predicted to be feasible by the surrogate models.
The more negative $\chmax(\s)$ is,
	the more likely is $\s$ to be feasible when evaluated by the blackboxes. 
Decreasing progressively the value of $\chmargin$ after each call of this 
	selection method will favor candidates that are increasingly likely to be feasible
	(but possibly with worse objective function values). 

\begin{myAlgo}[!ht]
	\algoNew
	{ 
		\algoLine{$\sbest \leftarrow \sinf$}

		\algoOpen{if first use of Method~\ref{algo:method4}, then:}
		{
			\algoLine{$\chmargin \leftarrow \min \{   0,   \underset{\substack{\s\in\Xh\\ \chmax(\s)<0}}{\max}\chmax(\s)    \} $}
		}
		\algoEnd 

		\algoOpen{for all $\s \in \Xh$, do:}
		{
			\algoOpen{if $\chmax(\s) \le \chmargin$ \AND $\fh(\s) < \fh(\sbest)$ \AND $d(\s,\XUS) > \Delta^\mesh$, then:}
			{
				\algoLine{$\sbest \leftarrow \s$}
			}
			\algoEnd
		}
		\algoEnd

		\algoOpen{if $\sbest \neq \sinf$, then:}
		{
			\algoLine{$\S \leftarrow \S \cup \{\sbest\}$}
			\algoLine{$\chmargin \leftarrow 2\,\chmax(\sbest)$}
		}
		\algoEnd
	}
	\caption{Selection of the best point with a constraint on the feasibility (Method~\ref{algo:method4})}
	\label{algo:method4}
\end{myAlgo}	

Note that this method requires that $\chmargin \le 0$ is always satisfied. Moreover, we also assume
that there is at least one $\s \in \Xh$ that is feasible.
If it is not the case,
which may happen in the first iteration,
	we will end up with $\chmargin > 0$
	and an inappropriate candidate will be selected.
To avoid this,
	we initialize $\chmargin$ to $0$ so that 
	the method may fail to return a candidate if that is the case.
	
	~\\


\newcommand{\iso}{\text{iso}}

{\bf Method~\ref{algo:method5}.}
The isolation distance is used here to detect local minima of the surrogate problem.
This concept is inspired from the \emph{topographic isolation} of a mountain summit,
	which measures the local significance of a summit.
It is defined as the distance to the closest higher summit.\footnote
	{
		In mountaineering,
			the topographic isolation of Mount Everest is infinite 
			and the summit with the second highest isolation is the Aconcagua in Argentina.
		The Aconcagua is not the second highest summit
			but there is no higher mountain in a 16,518 km range,
			making it the most important summit in the Americas
			and in the southern hemisphere.
	}

Transferred to optimization,
	the concept of topographic isolation is used to quantify the importance of a local minimum.
Its strict application is however impossible
	since it will require to prove that no other point within a certain distance of $\x$ is better than $\x$.
We can only compute isolation distance of the already evaluated points.
Consequently,
	we define the isolation distance as being 
	the distance from $\s$ to the closest point of $\Xh$ that is better than $\s$
	\begin{equation}
		d_{\iso}(\s) = \underset{\substack{\s'\in\Xh\\ \s'\!\prech\!\s}}{\min} \;  d\big( \s , \s' ).
	\end{equation}
Constraints are taken into account 
	by using the order relationship defined in Equation~\eqref{eq:order}. 
As a convention,
	if no point of $\Xh$ is better than $\s$,
	then $d_{\iso}(\s) = +\infty$.
With this definition, 
	the point of $\Xh$ with the highest isolation distance is also the best candidate in $\Xh$.
However,
	we have observed that
	the other points with a high isolation distance are often poor points far from any other point of $\Xh$.
To address this problem,
	we define the \emph{isolation number} of $\s \in \Xh$ as 
	the number of points of $\Xh$ within the ball of centre $\s$ and radius $d_{\iso}(\s)$
	\begin{equation}
 		n_{\iso}(\s) = \card \big\{ \s' : \s' \in \Xh, d(\s,\s') < d_{\iso}(\s) \big\}.
	\end{equation}
To have a high isolation number,
	a point must be better than many of its neighbors,
	which means that this criterion allows to detect local minima.
Note that Equation~\eqref{eq:sinf} implies that
	$d_{\iso}(\sinf) = n_{\iso}(\sinf) = 0$.
Method~\ref{algo:method5}  selects the point of $\Xh$
	that has the highest isolation number not yet in $\XUS$.
%

\begin{myAlgo}[!ht]
	\algoNew
	{ 
		\algoLine{$\sbest \leftarrow \sinf$}

		\algoOpen{for all $\s \in \Xh$, do:}
		{
			\algoOpen{if $n_{\iso}(\s) > n_{\iso}(\sbest)$  \AND  $d(\s,\XUS) > 0$, then:}
			{
				\algoLine{$\sbest \leftarrow \s$}
			}
			\algoEnd
		}
		\algoEnd

		\algoOpen{if $\sbest \neq \sinf$, then:}
		{
			\algoLine{$\S \leftarrow \S \cup \{\sbest\}$}
		}
		\algoEnd
	}
	\caption{Selection of the most isolated point (Method~\ref{algo:method5})}
	\label{algo:method5}
\end{myAlgo}


\newcommand{\density}{\text{density}}

{\bf Method~\ref{algo:method6}.}
The purpose of this method is to select points in neglected areas of the design space.
To do so,
	it selects points in areas heavily explored while solving~\eqref{eq:SurrogateProblem}
	but overlooked when solving~\eqref{eq:MainProblem}.
The \emph{density number} of $\s \in \Xh$ is defined as
	\begin{equation}
  		n_{\density}(\s) = \card \big\{ \s' : \s' \in \Xh, d(\s,\s') < d(\s,\XUS)  \big\}.
	\end{equation}
Method~\ref{algo:method6}  selects the point of $\Xh$ with the highest density number.
Note that, as for $n_{\iso}$, Equation~\eqref{eq:sinf} implies that $n_{\density}(\sinf)=0$.

\begin{myAlgo}[!ht]
	\algoNew
	{ 
		\algoLine{$\sbest \leftarrow \sinf$}

		\algoOpen{for all $\s \in \Xh$, do:}
		{
			\algoOpen{if $n_{\density}(\s) > n_{\density}(\sbest)$, then:}
			{
				\algoLine{$\sbest \leftarrow \s$}
			 }
			\algoEnd
		}
		\algoEnd

		\algoOpen{if $\sbest \neq \sinf$, then:}
		{
			\algoLine{$\S \leftarrow \S \cup \{\sbest\}$}
		}
		\algoEnd
	}
	\caption{Selection of a point in a populated area (Method~\ref{algo:method6})}
	\label{algo:method6}
\end{myAlgo}

\section{Parallel computing implementation}
\label{sec:algo}

We now describe how we implement the proposed SEARCH step.
We start with the surrogate models
	and follow with the algorithms used for solving \eqref{eq:SurrogateProblem}.
We conclude with how all of this is integrated with the MADS algorithm to solve \eqref{eq:MainProblem}.

\subsection{Surrogate models}

Several surrogate models are mentioned in~\cite{Haftka2016} regarding their use in parallel computing approaches, including
	Kriging,
	radial basis functions (RBF),
	support vector regression (SVR),
	and polynomial response surfaces (PRSs).
For a more exhaustive description and comparison of surrogate models,
	see~\cite{AlizadehAllenMistree2020}.
Based on our previous work reported in \cite{Lowess2016},
	we choose to use the locally weighted scatterplot smoothing (LOWESS) surrogate modeling approach 
	\cite{Cle1979, Cle1981, Cle1988a, Cle1988b}.

LOWESS models generalize PRSs and kernel smoothing (KS) models.
PRSs are good for small problems,
	but their efficacy decreases for larger, highly nonlinear, or discrete problems. 
KS models tend to overestimate low function values
	and underestimate high ones,
	but usually predict correctly which of two points yields the smallest function value~\cite{Ensemble2016}.
LOWESS models build a linear regression of kernel functions around the point to estimate.
They have been shown
	to be suitable for surrogate-based optimization \cite{Lowess2016}.
Their parameters are chosen using an error metric called 
	\quote{aggregate order error with cross-validation} (AOECV),
	which favors equivalence between the original problem~\eqref{eq:MainProblem}
	and the surrogate problem~\eqref{eq:SurrogateProblem}~\cite{Ensemble2016}.
%
	We use the SGTELIB implementation of the surrogate models,
	which is now integrated as a surrogate library in NOMAD version 3.8~\cite{sgtelib}.
Specifically,
	the considered LOWESS model is defined in SGTELIB as follows.
	\begin{equation*}
 		\text{\tt TYPE Lowess DEGREE 1 RIDGE 0 SHAPE\_COEF OPTIM KERNEL\_TYPE OPTIM},
	\end{equation*}
	which means that the local regression is linear,
	the ridge (regularization) coefficient is 0, 
	and the kernel shape and kernel type are optimized to minimize the aggregate order error (see~\cite{Lowess2016}).

The LOWESS model is built as described in Appendix~\ref{sec:Lowess}.
Only the Gaussian kernel was considered in~\cite{Lowess2016}.
Six additional kernel functions have meanwhile been implemented in SGTELIB.
Accordingly,
	not only $\lambda$ (kernel shape) is chosen to minimize AOECV,
	but also $\phi$ (kernel type). 
%

\subsection{Surrogate problem solution}

The surrogate problem \eqref{eq:SurrogateProblem} is solved by means of an inner instance of MADS; it is
	initialized by a Latin hypercube search (LHS)~\cite{McCoBe79a,SaWiNoLH2003} and 
	uses variable neighborhood search (VNS)~\cite{MlHa97a,HaMl01a} as the SEARCH step
	and a large budget of function evaluations
	(10,000).
The POLL step is performed using the ORTHO~2N directions option.
This inner MADS  is implemented in the SEARCH step of the outer MADS.
%

The LHS guarantees that there are surrogate cache points widely spread over the design space.
To ensure this,
	30\% of all function evaluations (i.e., 3,000) are devoted to the LHS.
Four additional points are considered:
	The current best feasible point of the original problem,  
	the current best infeasible point of the original problem, 
	the best feasible point obtained by solving the most recent surrogate problem instantiation,
	and the best infeasible point by solving the most recent surrogate problem instantiation. 
These points are used as initial guesses of the inner MADS problem,
	which will be run until the remaining evaluation budget is exhausted.
This  budget will be shared between the POLL step and the VNS
	in a default proportion where 75\% is devoted to VNS, which favors the exploration of multiple local attraction basins.
A large number of evaluations that build the surrogate cache $\Xh$
	favors an accurate solution of the surrogate problem.
Using LHS, VNS, and a large number of function evaluations
	 ensures that $\Xh$ contains highly promising candidates
	for the solution of the original problem~\eqref{eq:MainProblem}.

\subsection{The modified MADS algorithm}
\label{sec:outerMADS}

Recall that each iteration of 
	MADS includes a SEARCH step (performed first) and a POLL step.
Let $t$ denote the current iteration.
Then,
	$\X_t$, $\Xh_t$, and $\S_t$ denote the sets $\X$, $\Xh$, and $\S$
	considered by MADS at iteration~$t$.
Let $q$ be the number of available CPUs,
	i.e., the number of blackbox evaluations that can be performed in parallel.
The proposed MADS for exploiting $q$ CPUs is as follows.
First, 
	the SEARCH step proceeds by solving the surrogate problem 
	to populate the set  $\Xh_t$.
From that set,
	$q$ candidates are selected and returned to be evaluated by the blackbox(es) in parallel.
The selection is made
	by cycling through a user-defined subset of the six proposed selection methods
	(Section~\ref{sec:detailedDescription})
	until a total of $q$ candidates are selected,
	or until all selection methods consecutively failed to add a candidate to $\S_t$.
If $q$ is smaller than the number of selection methods retained by the user,
	we do not necessarily go through all the  methods,
	but stop as soon as  we get $q$ candidates.
		
If/when the SEARCH step fails to return a better objective function value, the MADS algorithm proceeds to the POLL step. 
Let $\P_t$ be the set of candidates produced by the polling directions at the iteration $t$.
The cardinality of  $\P_t$  is denoted by  $|\P_t|$.
To be consistent with our need to fulfill continuously all the available CPUs with evaluations,
	additional candidates are added
	so that  $|\P_t|$  is at least $q$ or a multiple of $q$.
This is accomplished by means of NOMAD's intensification mechanism ORTHO~1.
If $|\P_t| = q$,
	all poll candidates of $\P_t$ are evaluated concurrently, eliminating the need to order them.
If $|\P_t| > q$,
	then the points in $\P_t$ are regrouped in several blocks of $q$ candidates.
The blocks are then evaluated sequentially and opportunistically,
	which means that if a block evaluation leads to a success,
	the remaining blocks are not evaluated.
To increase the probability of success from the first  block,
	and hence avoiding to proceed with the remaining ones,
	the candidates of $\P_t$ are sorted using the surrogate models 
	and distributed in the  blocks so that the more promising ones are in the first block.
%


Recall that $\mesh_t$ is the current mesh,
	$\Delta^\mesh_t$ is the associated mesh size parameter,
	$\Delta^{\text{P}}_t$ is the corresponding mesh poll parameter,
	and $\x^*_t$ is the best solution found at iteration $t$.
Finally,
	the set $\X_t$ is updated with all the points $\x$ in $\S_{t-1}$ and $\P_{t-1}$
	that have been evaluated during the previous iteration.
The process is summarized in Algorithm~\ref{algo:mads}.

\begin{myAlgo}[h]
\algoNew{ 
	\algoSection{[1] Initialization}{
		\algoLine{$t \leftarrow 0$}
		\algoLine{Set initial poll and mesh sizes $\Delta^{\text{P}}_0 \ge \Delta^\mesh_0 > 0$}
		\algoLine{Initialize $\X_0$ with starting points}
		\algoLine{Evaluate $\{f(\x),c_1(\x),c_2(\x),\hdots,c_m(\x)\} \; \forall \x \in \X_0$}
	}
				
	\algoSection{[2] Model search}{
		\algoLine{Use $\X_t$ to build $\fh$ and $ \{ \ch_j \}_{j \in J}$}
		\algoLine{Solve surrogate problem~\eqref{eq:SurrogateProblem} using the inner MADS instance}      
		\algoLine{$\Xh_t \leftarrow $ Set of points evaluated with surrogate model while solving~\eqref{eq:SurrogateProblem}}
		\algoLine{$\S_t \leftarrow $ Cycle through selection steps to select $q$ points of $\Xh_t$}
		\algoLine{$\S_t \leftarrow $ Projection of the points of $\S_t$ onto mesh $\mesh_t$}
		\algoLine{Parallel evaluation of $\{f(\x),c_1(\x),c_2(\x),\hdots,c_m(\x)\} \; \forall \x \in \S_t$}
		\algoLine{If success, \algoFont{goto} \textbf{[4]}}
	}

	\algoSection{[3] Poll}{
 		\algoLine{Build poll set $\P_t$}
		\algoLine{Sort $\P_t$ according to $\fh$ and $ \{ \ch_j \}_{j \in J}$}
		\algoLine{Parallel evaluation of $\{f(\x),c_1(\x),c_2(\x),\hdots,c_m(\x)\} \; \forall \x \in \P_t$}
	}
	
	\algoSection{[4] Updates}{
		\algoLine{$t\leftarrow t+1$}
		\algoLine{Update $\Delta^\mesh_t$, $\Delta^{\text{P}}_t$, $\x^*_t$ and $\X_t$}
		\algoLine{If no stopping condition is met, \algoFont{goto} \textbf{[2]}}
	}
}
\caption{The proposed \mads optimization algorithm}
\label{algo:mads}
\end{myAlgo}

\section{Numerical investigation}
\label{sec:test}

The proposed SEARCH step technique is tested using five optimization problems.
We first describe the algorithms considered for benchmarking.
Next, numerical results are presented and discussed for
	the five engineering design problems. 

\subsection{Compared algorithms}

Five solvers are compared in our numerical experiments, all based on the MADS algorithm
	and implemented using NOMAD 3.8~\cite{Le09b}.
	This ensures avoiding coding biases
	since features are identical among solvers.
	\begin{itemize}
		\item	{\bf MADS.}
			Refers to the POLL step of MADS,
				without any SEARCH step,
				where $2n$ directions are generated and evaluated in parallel.
			If needed,
				$k$ additional directions are generated such that $2n+k$ is a multiple of $q$.

		\item	{\bf Multi-Start.}
			Consists of $q$ parallel  runs of MADS.
			They are totally independent and each instance runs on its own CPU.
			Each instance proceeds to its evaluations sequentially,
				one after the other.
			Only the POLL step is executed and
				no cache is shared between running instances.

		\item	{\bf LH Search.}
			The MADS solver mentioned above
				using a Latin hypercube search (LHS) at the SEARCH step,
				where $q$ candidates are generated 
				and evaluated in parallel.

		\item	{\bf Lowess-A.}
			The MADS solver mentioned above
				with the described surrogate optimization conducted at the SEARCH step.
			The $q$ candidates are selected by cycling through
				Methods~\ref{algo:method1} and \ref{algo:method2},
				and then evaluated in parallel.

		\item	{\bf Lowess-B.}
			The MADS solver mentioned above
				with the proposed surrogate optimization conducted at the SEARCH step.
			The $q$ search candidates are selected by cycling through
				Methods~\ref{algo:method3}, \ref{algo:method4}, \ref{algo:method5}, and \ref{algo:method6},
				and then evaluated in parallel.
	\end{itemize}

Both LOWESS solvers are exactly like Algorithm~\ref{algo:mads},
	excepted for the used selection methods.
The only difference between them and the LHS solver 
	is that the surrogate optimization approach is replaced by a LHS at the SEARCH step.
This should allow us to determine whether surrogate optimization has any advantage over a random search.
The MADS solver is used as the baseline.
Finally,
	the Multi-Start solver is considered to ensure 
	that one should not proceed with $q$ independent narrow trajectories
	instead of one single trajectory having $q$ wide evaluations.

\subsection{Engineering design optimization problems}

The above solvers are compared on five engineering design application problems.
A short description follows below for each problem.
More details are provided in Appendix~\ref{annexe:detailedProb}.

	\begin{itemize}
		\item	{\bf TCSD.}
			The Tension/Compression Spring Design problem
				consists of minimizing the weight of a spring under mechanical constraints~\cite{Garg2014, Arora2004, Belegundu1982}.
			This problem has three variables and four constraints.
			The design variables define the geometry of the spring.
			The constraints concern  shear stress, surge frequency,  and minimum deflection.
	
		\item	{\bf Vessel.}
			This  problem
				considers the design of a compressed air storage tank
				and has four design variables and four constraints~\cite{Garg2014, KaKr1993}.
			The  variables define the geometry of the tank
				and the constraints are related to the volume,
				pressure,
				and solidity of the tank.
			The objective is to minimize the total cost of the tank,
				including material and labour. 
	
		\item	{\bf Welded.}
			The welded beam design problem (Version I)
				has four variables and six constraints~\cite{Garg2014, rao1996}.
			It aims at minimizing the construction cost of a beam,
				under shear stress,
				bending stress,
				and deflection constraints.
			The design variables define the geometry
				and the characteristics of the welded joint.
	
		\item	{\bf Solar~1.}
			This optimization problem aims at maximizing the energy received over a period of 24 hours
				under five constraints related to budget
				and heliostat field area~\cite{MScMLG}.
			It has nine variables, including 
an integer one without an upper bound.
	
		\item	{\bf Solar~7.}
			This  problem aims at maximizing the efficiency of the receiver over a period of 24 hours
				for a given heliostats field
				under six binary constraints~\cite{MScMLG}.
			It has seven design variables, including 
an integer one without an upper bound.
	\end{itemize}

A progressive barrier is used to deal with the aggregated constraints~\cite{AuDe09a}.
The three first problems are easier relative to the last two ones.
However,	it is difficult to find a feasible solution for the TCSD problem.
Among all the considered problems,
	Solar~1 is certainly the most difficult one.

\subsection{Numerical experiments}

We compare the efficiency of each  solver for different values of block size
	$q \in \{1, 2, 4, 8, 16, 32, 64\}$.
As an example,
	we will use \quote{Lowess-A~16} to refer to the solver that relies on LOWESS models
	cycling over Methods~\ref{algo:method1} and~\ref{algo:method2}
	considering a block size $q=16$.
For each problem,
	we generated 50 sets of 64 starting points with Latin hypercube sampling~\cite{McCoBe79a}.
For all solvers other than \quote{Multi-Start},
	only the first point of each set is used to perform  optimizations.
Doing so,
	we get 50 runs from the same starting points for each solver,
	 each problem,
	and each value of $q$.
For \quote{Multi-Start},
	since $q$ independent and parallel sequential runs of MADS must be performed,
	we use the $q$ first points of each set.
Doing so,
	we still get 50 runs for each problem and each $q$,
	while ensuring that all starting points are the same for all solvers.

To avoid that all \quote{LH Search} runs end up with nearly identical solutions for a given $q$,
	we use a random seed for initializing each LHS.
%

For the relatively  three simpler problems (TCSD, Vessel, and Welded),
	a budget of 100 block evaluations is allocated.
For the two relatively difficult problems (Solar~1 and Solar~7),
	the budget is increased to 200 block evaluations.
This means that,
	for a given problem,
	all solvers will have the same \quote{wall-clock} time,
	but  not necessarily the same resources
	(number of CPUs available for block evaluations)
	nor the same total number of blackbox evaluations.


\subsubsection*{Solution quality}

Figures~\ref{fig:summary1} and~\ref{fig:summary2}
	represent the distribution of the final objective function over the 50 runs
	for each problem, each solver, and each   block size $q$.
The minimum and maximum objective values that we obtained from the runs are indicated by circles in the figures.
Lower and  upper  quantiles are delimited by boxes.
Median values are represented by a bar into the boxes.
The more a distribution is  on the left side, the
	better the combination of solver and $q$ is.
Since we are mostly interested in the best combinations,
	the figures only focus on the smallest values.
Otherwise,
	it would be difficult to discern the difference among the best combinations.
All combinations for which the distribution is cut on the right side are performing poorly.

\renewcommand{\figWidth}{0.99\linewidth}
\begin{figure}[p]  
 	\center
	\includegraphics[width=\figWidth]{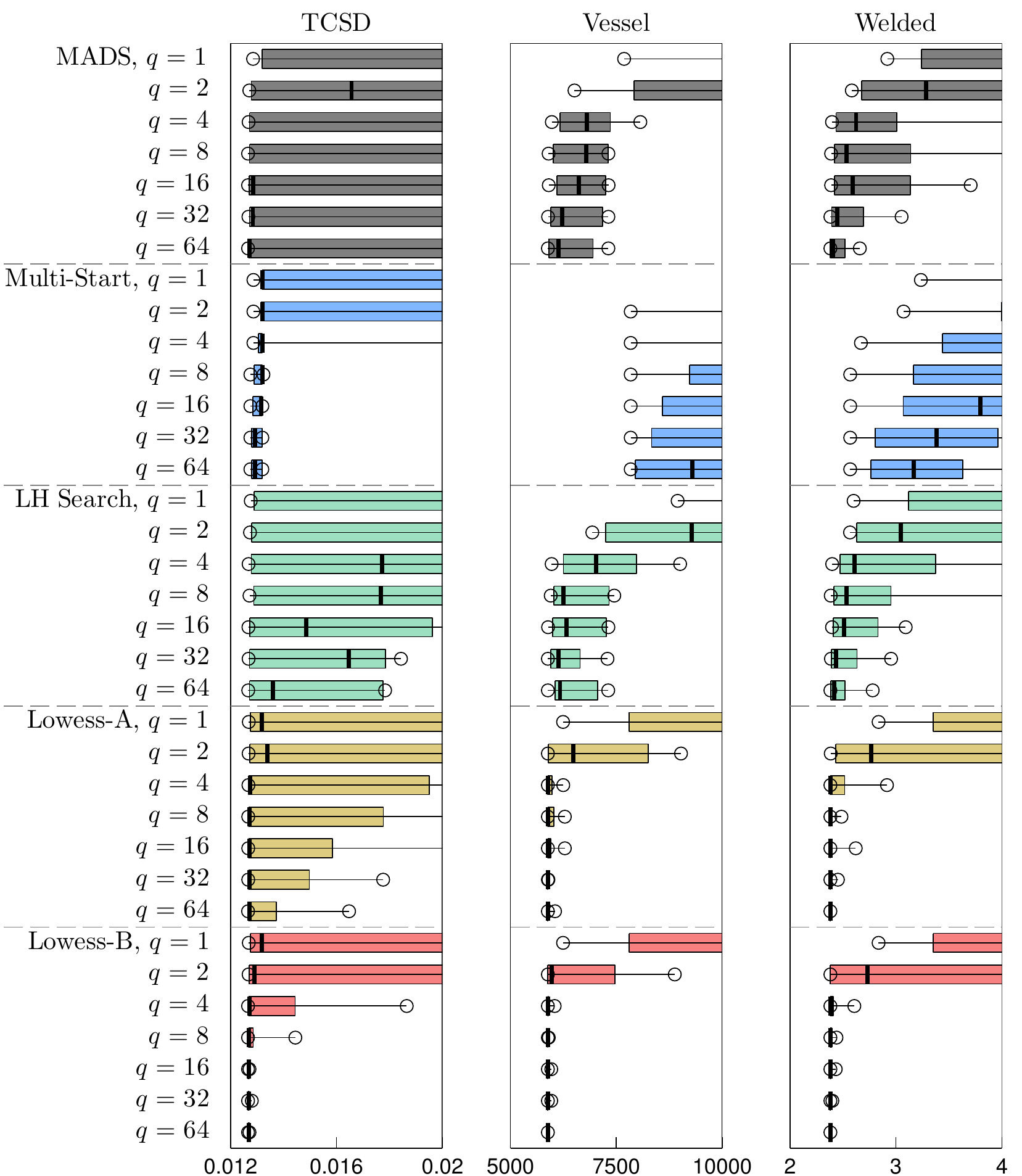}
	\caption{Performance summary for the TCSD, Vessel, and Welded problems over 50 runs}
	\label{fig:summary1}
\end{figure}

\renewcommand{\figWidth}{0.99\linewidth}
\begin{figure}[p]  
	\center
	\includegraphics[width=\figWidth]{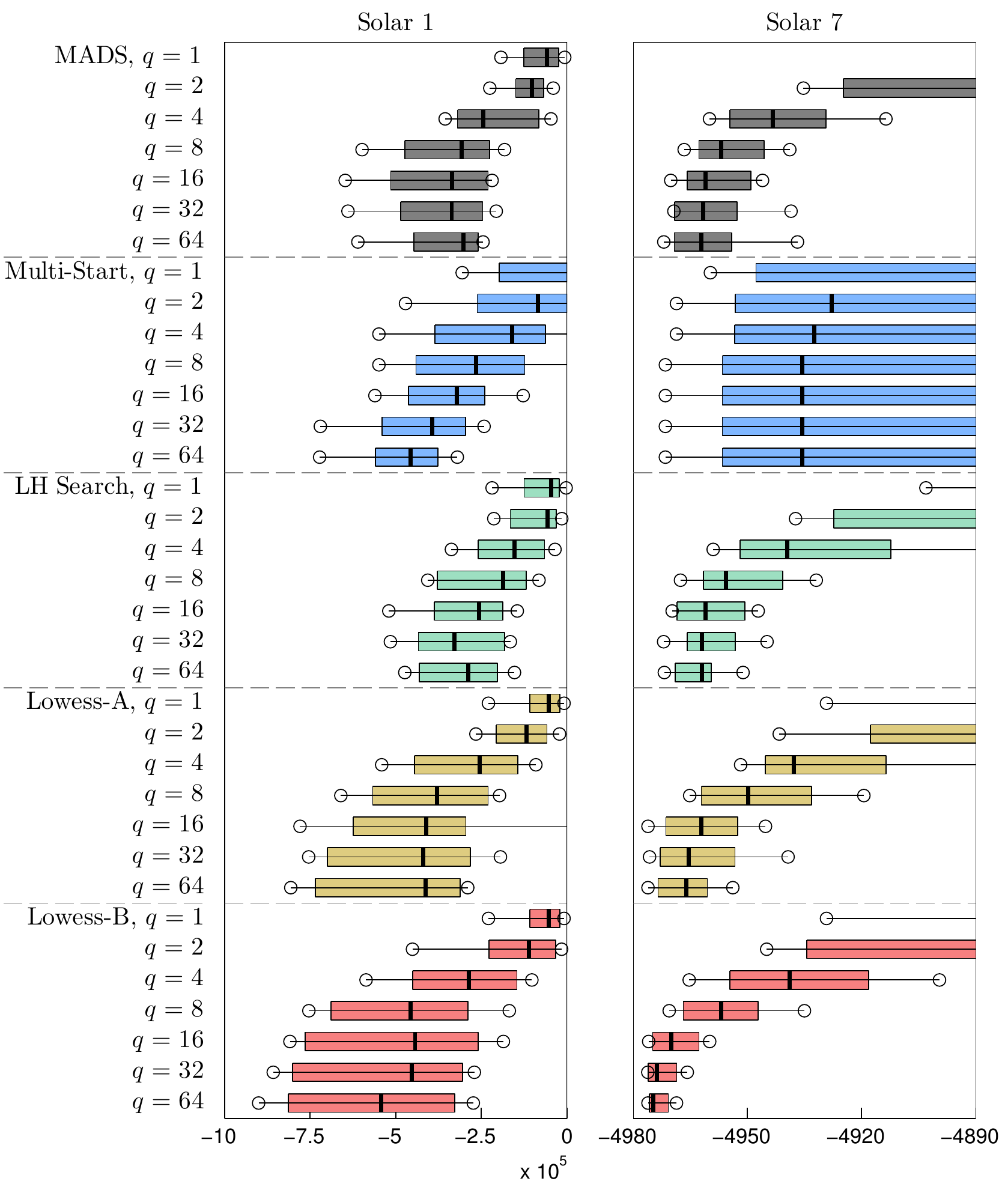}
	\caption{Performance summary for the Solar~1 and Solar~7 problems over 50 runs}
	\label{fig:summary2}
\end{figure}

For the three simpler problems
	(TCSD, Vessel, and Welded, Figure~\ref{fig:summary1}),
	the LOWESS solvers
	(and in particular \quote{Lowess-B})
	are by far superior to the solvers that do not rely on surrogate optimization.
	
For the TCSD problem,
	the \quote{MADS} solver often failed to find a feasible design
	(thus leading to infinite objective values),
	even with a large number of evaluations per block.
The four other solvers always managed to find a feasible point for at least 75\% of the runs. The \quote{Lowess-B} solver performs better than any of the other ones.
We see that \quote{Lowess-A~64} is outperformed by \quote{Lowess-B~8}. 
As the TCSD problem is very constrained,
	the final objective function value depends on the initial guess.
This is why the \quote{Multi-Start} solver performs quite well on this problem.

The same trend is observed for the Vessel, and Welded problems.
	\quote{Lowess-B} performs better than \quote{Lowess-A}, which outperforms \quote{LHS} or \quote{MADS}.
In particular, 
	\quote{Lowess-B 8} outperforms the solvers \quote{Lowess-A 8/16/32}.
As expected,
	increasing the block size improves performance.
However,
	for the \quote{Lowess-B} solver
	these three problems are easy to solve, so it is difficult to see an advantage of using parallel computing
	because the global optimum is found most of the time
	within 100 block evaluations for a block size of 16 or more.
The \quote{Multi-Start} solver performs rather poorly on these two problems.

The numerical results generally follow the same trend for the two Solar problems (Figure~\ref{fig:summary2}).
For a block of equal size,
	the LOWESS solvers outperform the other solvers
	while \quote{Lowess-B} outperforms \quote{Lowess\nobreakdash-A}. 

\subsubsection*{Convergence rate}

We now examine the convergence rate of the solvers for the case where $q = 64$.
Figures~\ref{fig:curves-median-1} and~\ref{fig:curves-median-2}
	depict the evolution of the median objective function value of 50 runs
	as a function of the number of block evaluations.
For each problem,
	the plots on the left compare the convergence of the five solvers
	with  blocks of size  $q = 64$
while the plots on the right compare the convergence of the best-performing solver i.e., \quote{Lowess-B}, for block sizes ranging from $q = 1$ to 64.
%

\renewcommand{\figWidth}{0.49\linewidth}
\begin{figure}[ht!]
	\center
		TCSD problem\linebreak
		\includegraphics[width=\figWidth]{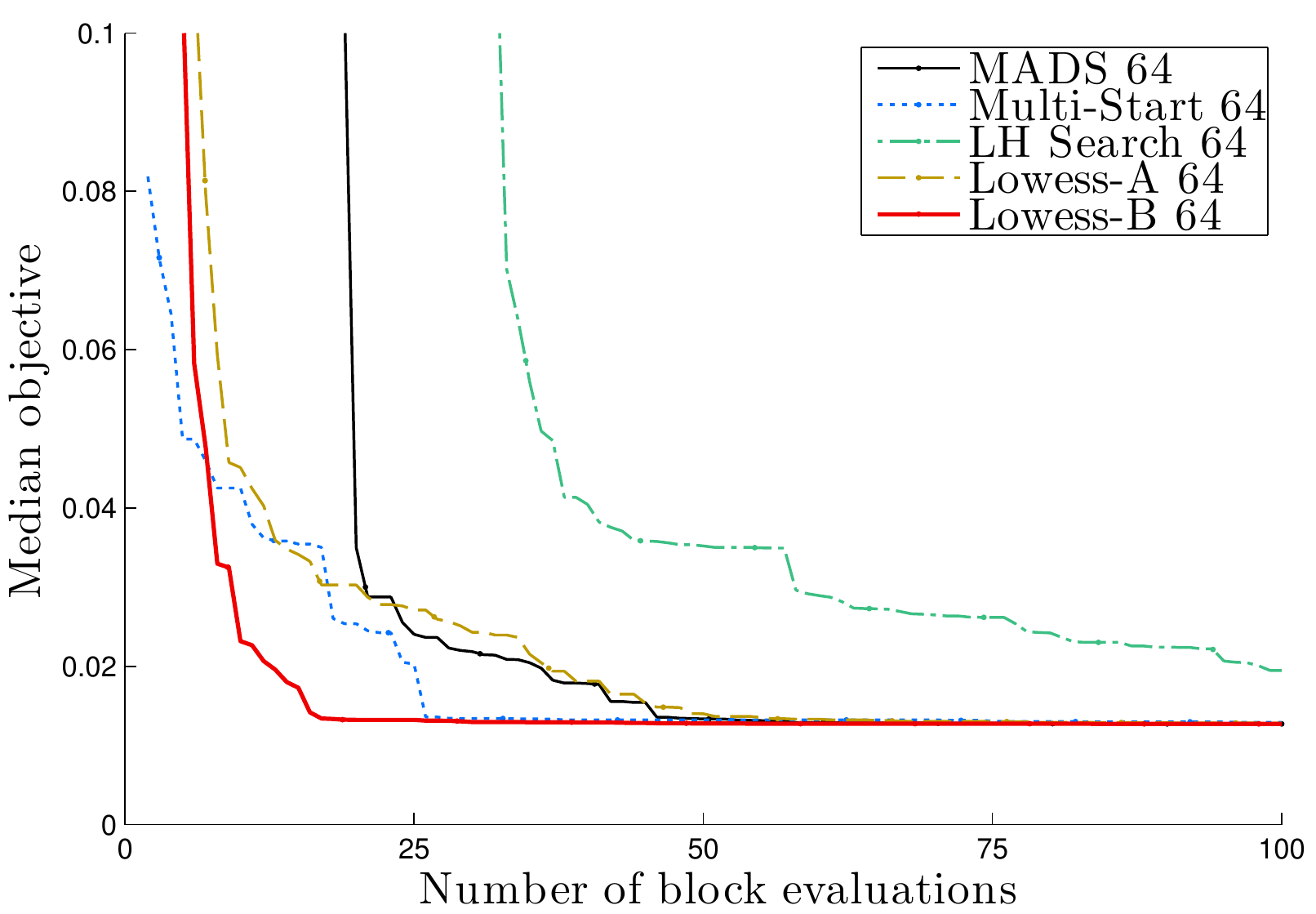}
		\includegraphics[width=\figWidth]{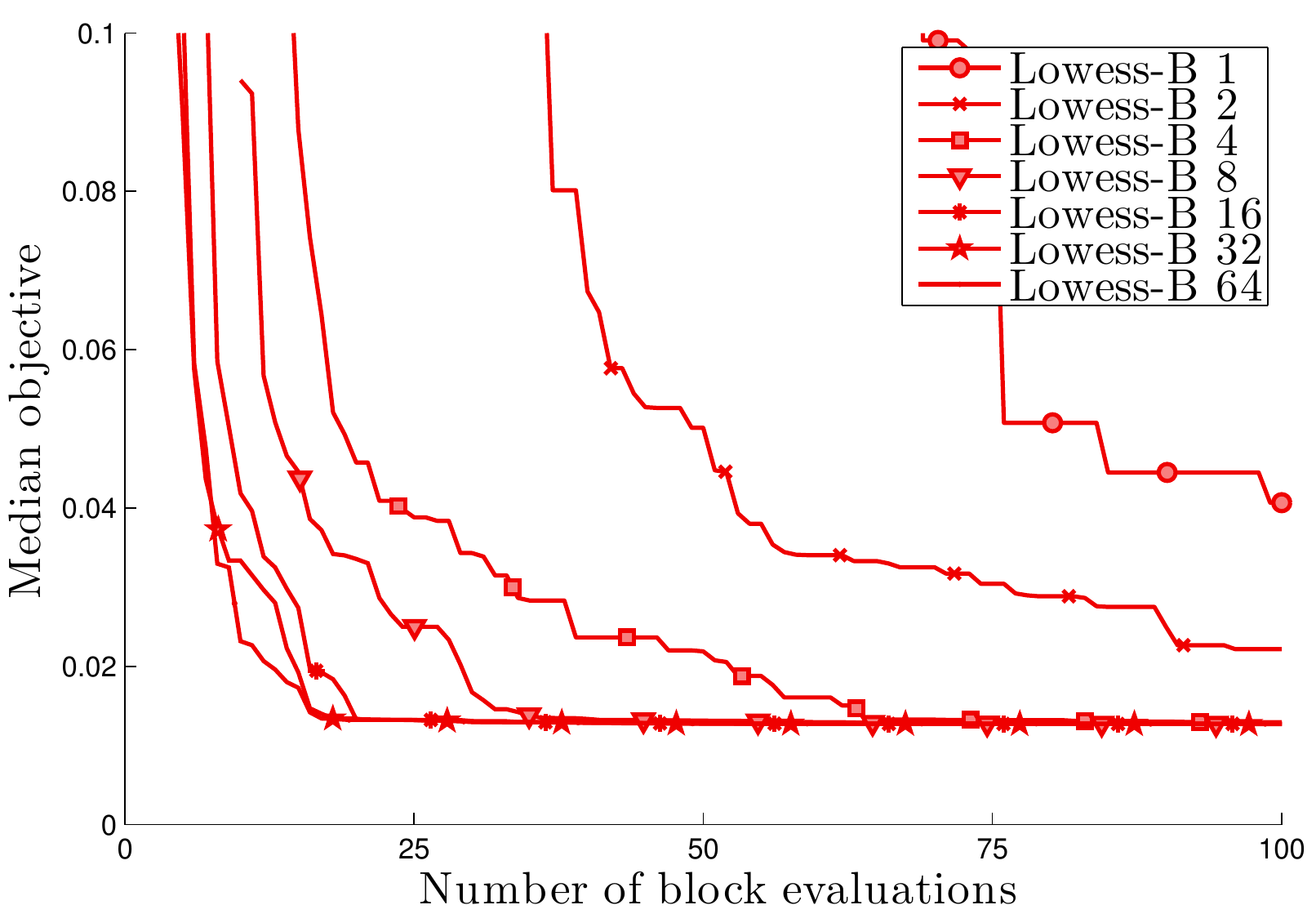}\newline

		Vessel problem\linebreak
		\includegraphics[width=\figWidth]{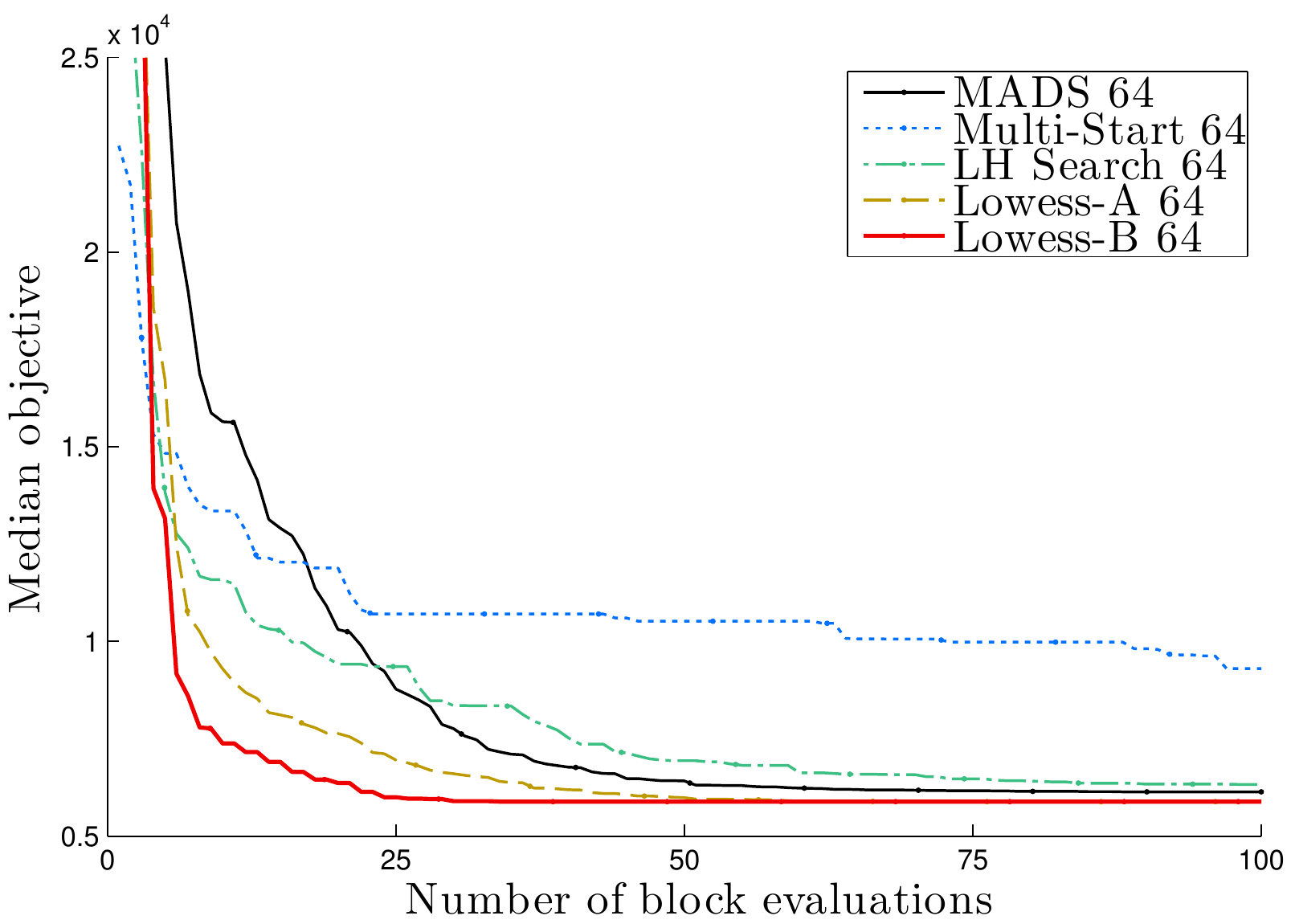}
		\includegraphics[width=\figWidth]{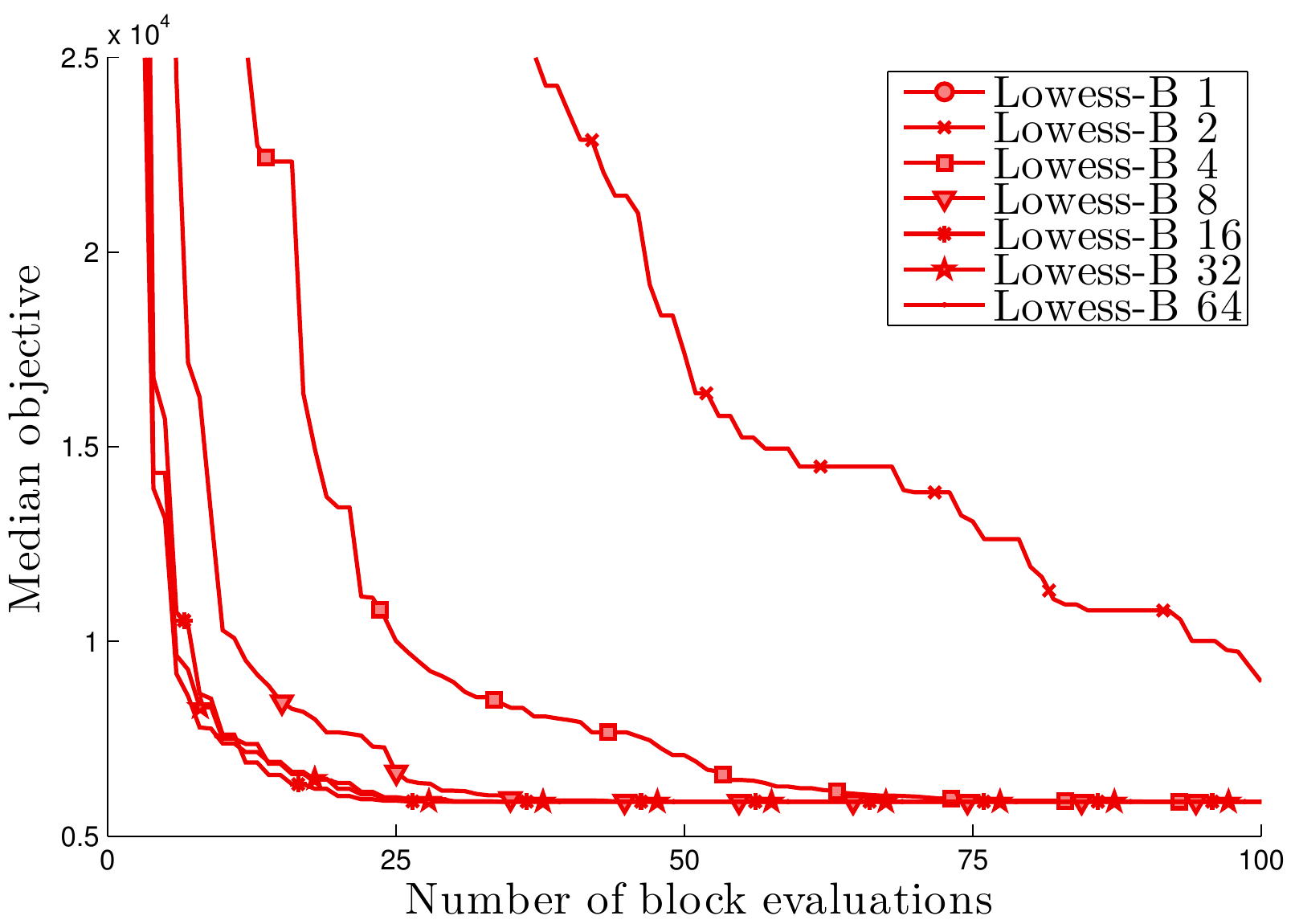}\newline

		 Welded problem\linebreak
		\includegraphics[width=\figWidth]{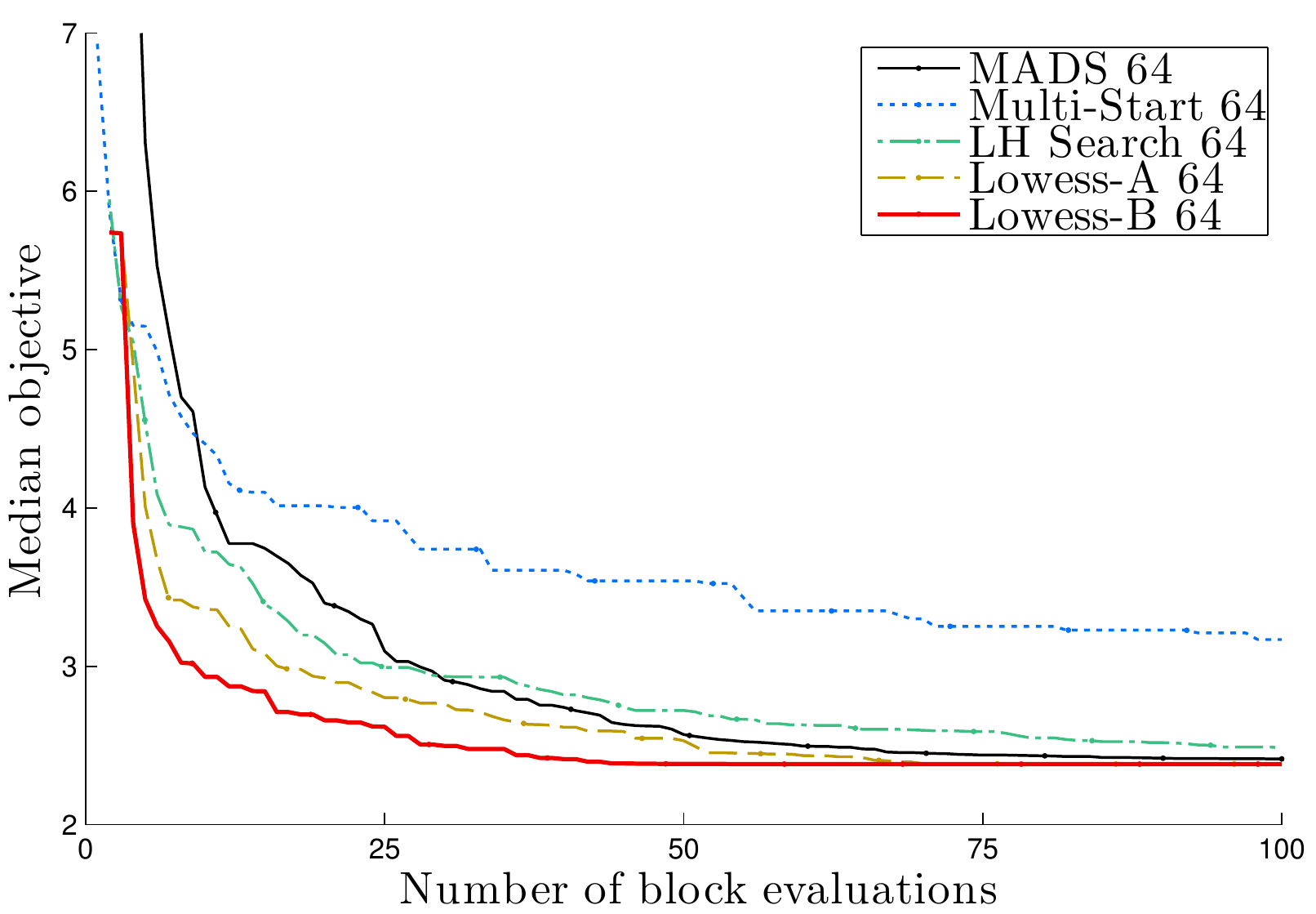}
		\includegraphics[width=\figWidth]{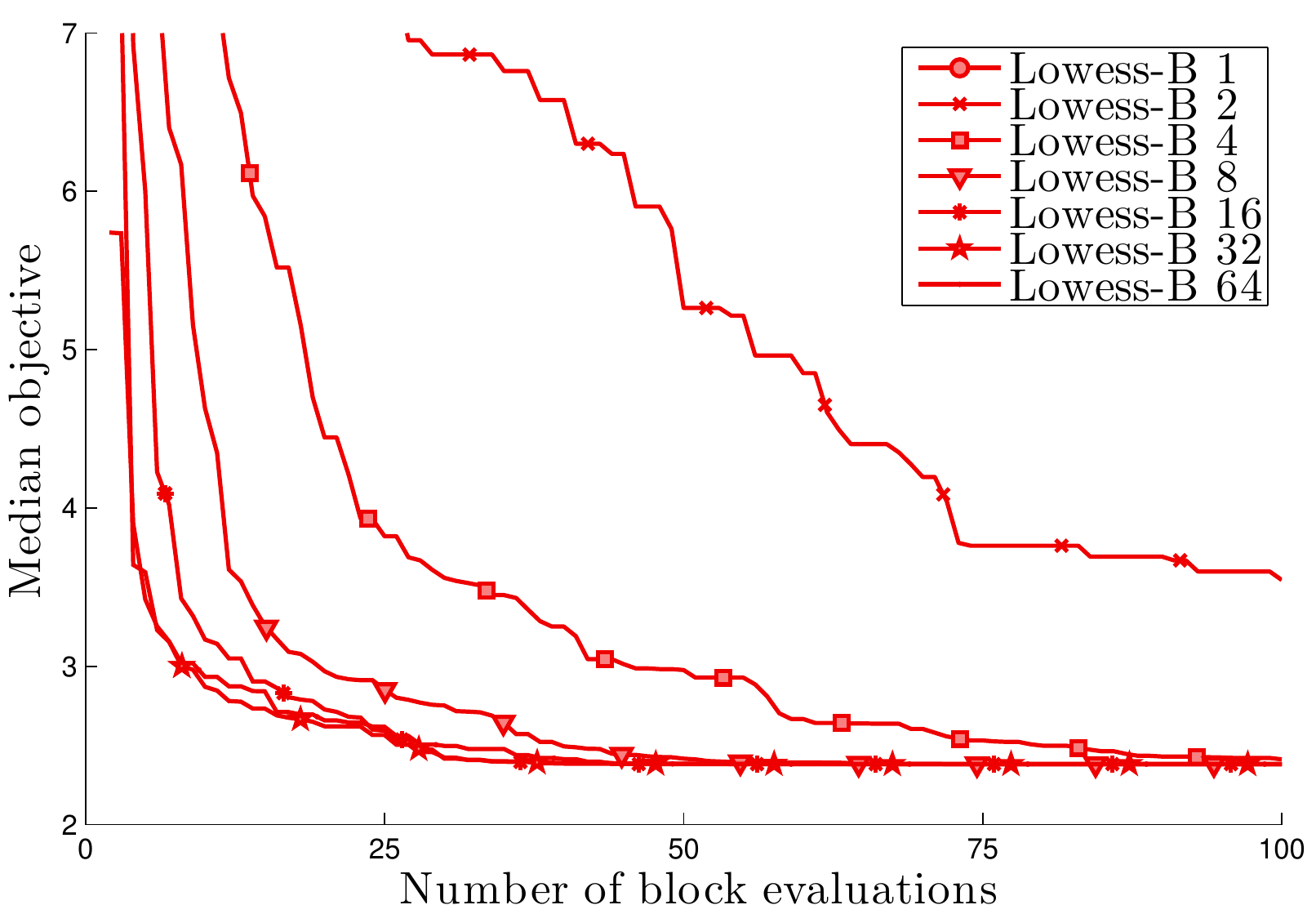}\newline
	\caption{Results for the TCSD, Vessel, and Welded problems; median objective value of 50 runs }
	\label{fig:curves-median-1}
\end{figure}

\renewcommand{\figWidth}{0.49\linewidth}
\begin{figure}[ht!]
	\center
		Solar~1 \linebreak
		\includegraphics[width=\figWidth]{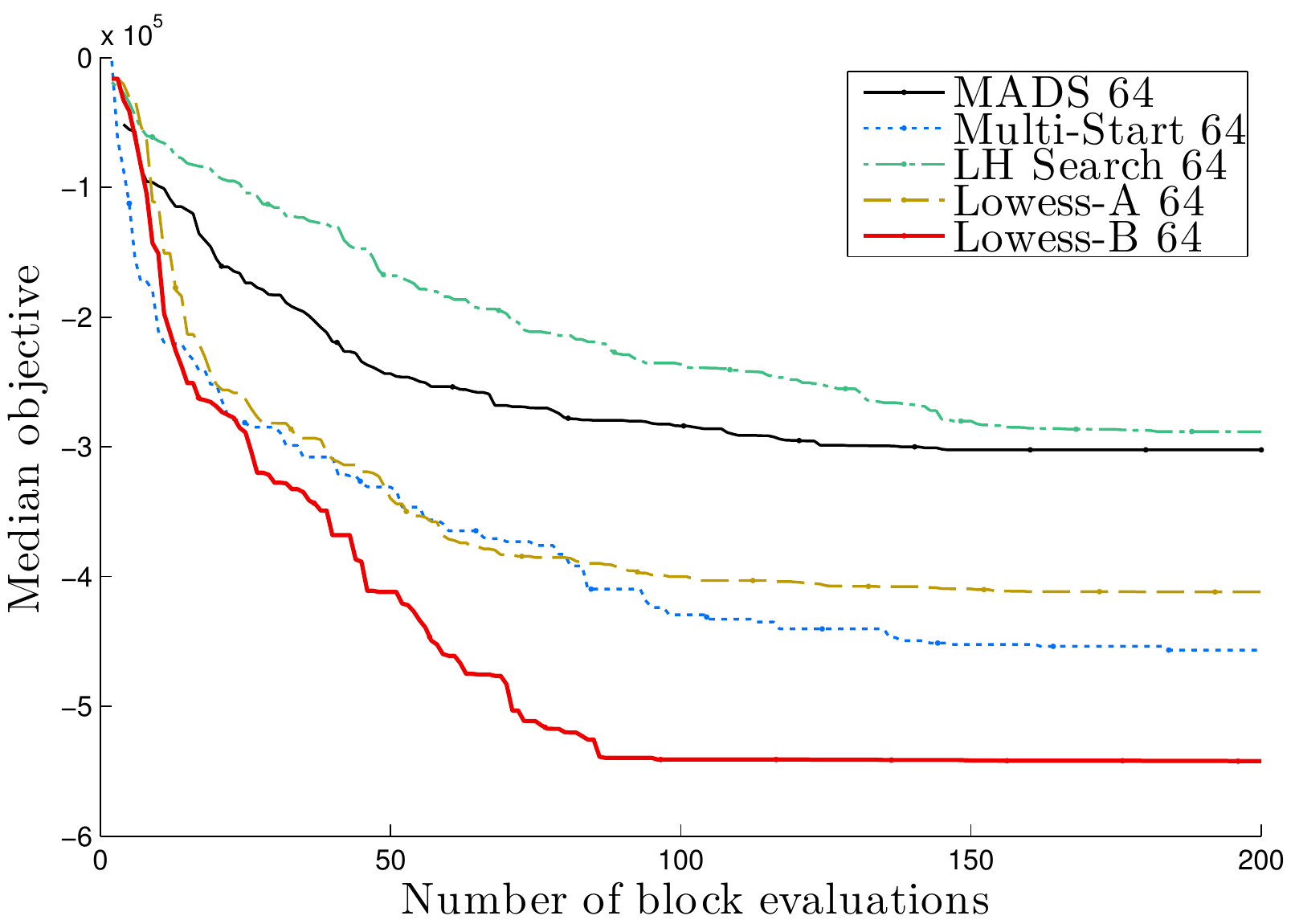}
		\includegraphics[width=\figWidth]{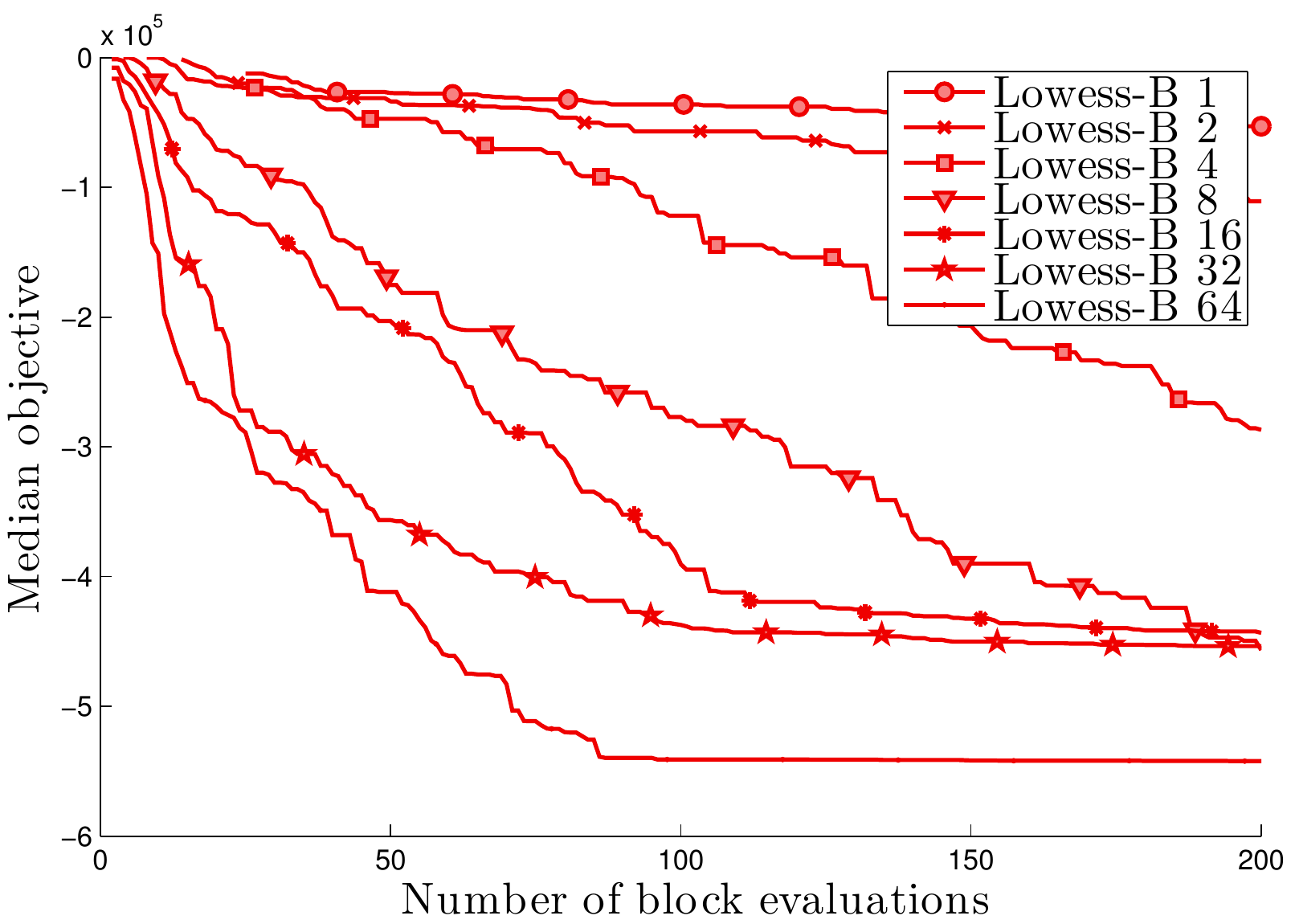}\newline

		Solar~7 \linebreak
		\includegraphics[width=\figWidth]{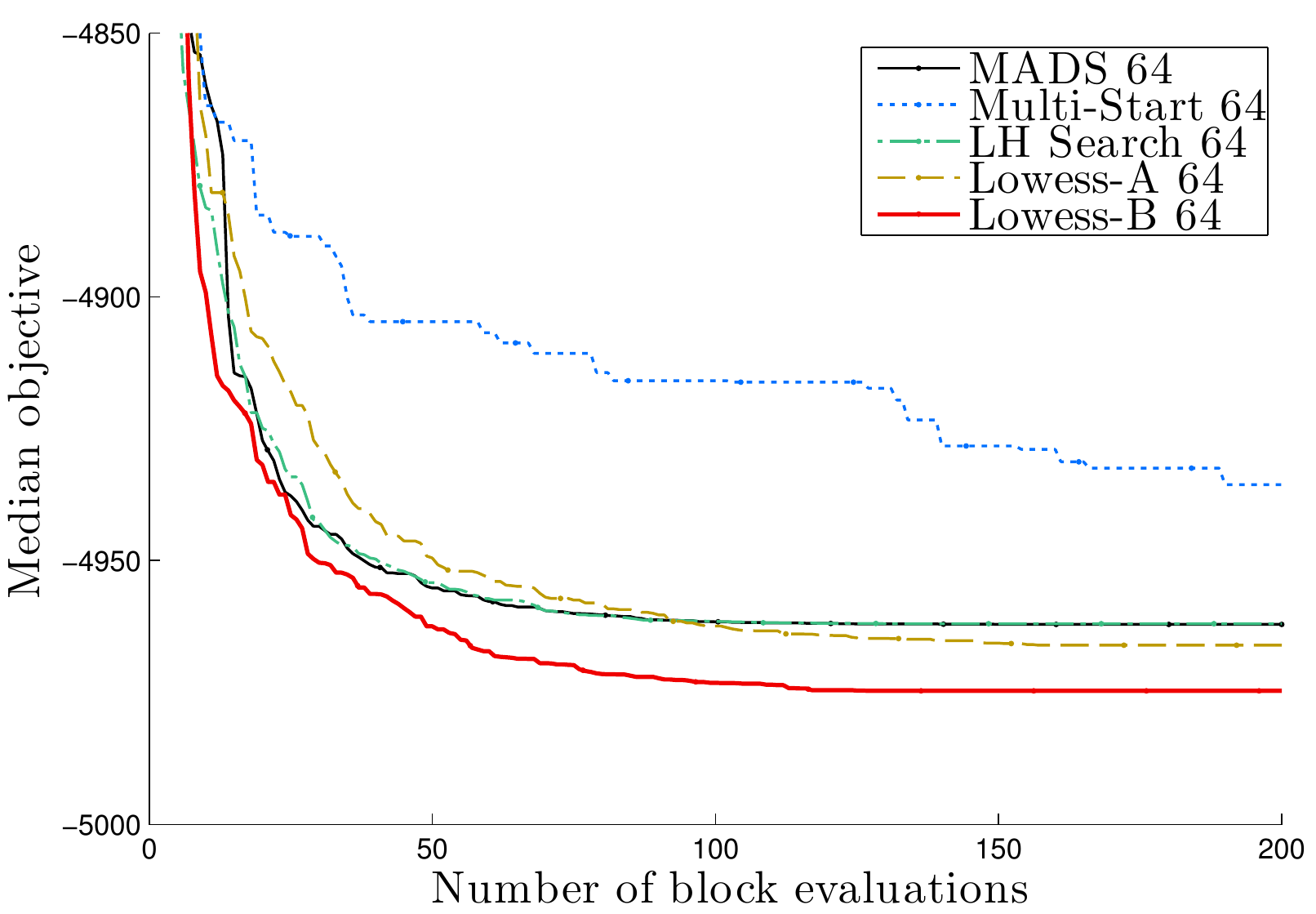}
		\includegraphics[width=\figWidth]{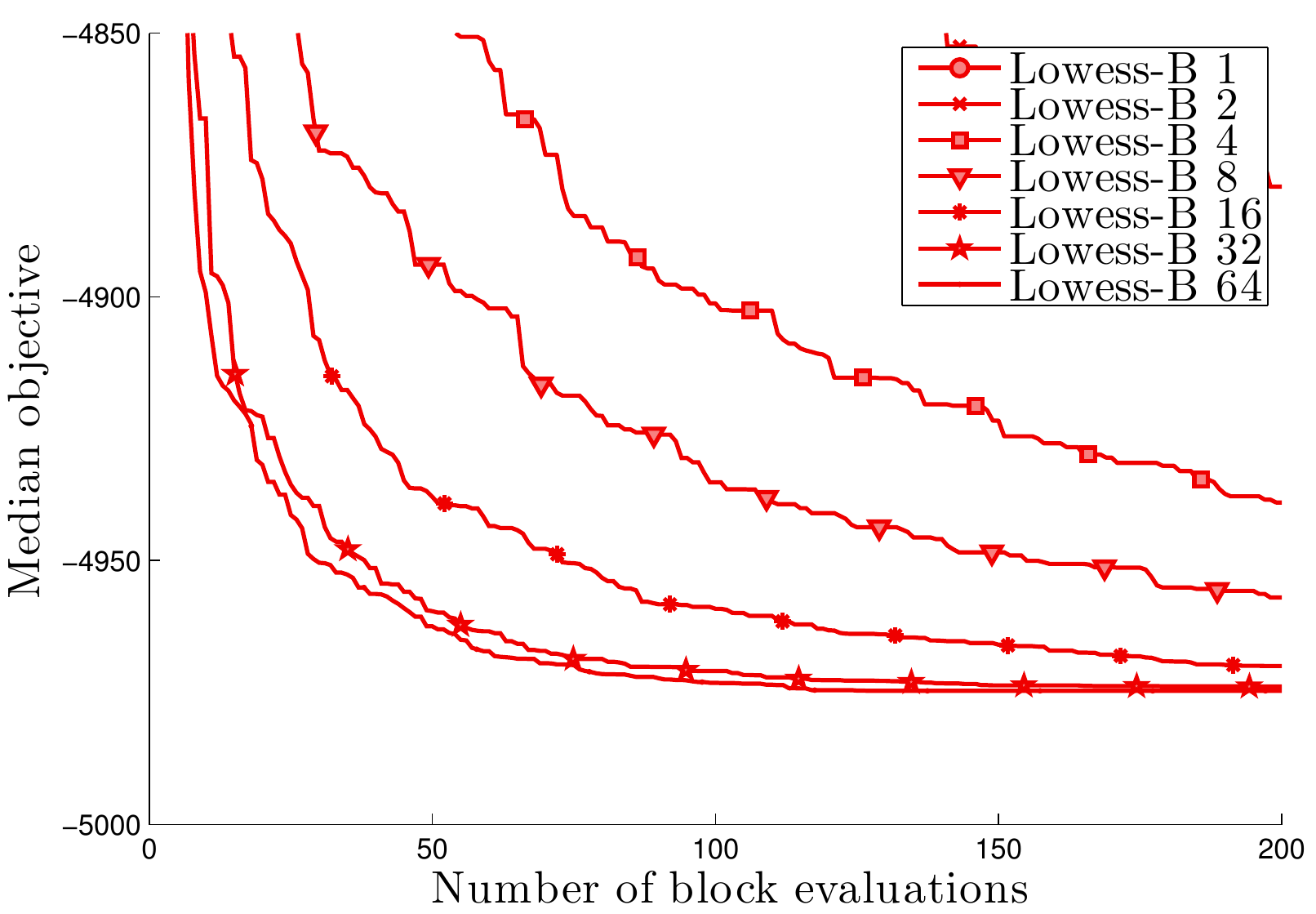}\newline
	\caption{Results for the Solar~1 and Solar~7 problems; median objective value of 50 runs}
	\label{fig:curves-median-2}
\end{figure}

We can conclude that \quote{Lowess-B} yields the best solutions
faster than any other solver
	(for $q = 64$).
The worst-performing solvers are \quote{Multi-Start} for problems TCSD and Solar~1
	and \quote{LH Search} for problems Vessel, Welded, and Solar~7.
It is also notable that although \quote{MADS} does not use the SEARCH step, 
	it performs generally well, 
	except for Solar~1.
\quote{Lowess-A} performed well but does not clearly outperform other solvers.
%

Considering the performance of \quote{Lowess-B} as a function of $q$, we observe that, as expected, 
convergence improves for larger values of $q$.
Depending on the problem,
	there may be a saturation point beyond which an increase of $q$ does not effect an improvement.
E.g., a saturation point arises around $q = 8$ or 16 for TCSD, Vessel, and Welded.
On the contrary,  $q$ could be even larger than 64 for Solar~1 as more CPUs can be utilized.
%
%

\subsubsection*{Performance profiles}

We now consider performance profiles, which indicate
	the percentage of runs where the problem is solved within a deviation from the best known solution $\tau$ under a budget of function evaluations ~\cite{DoMo02}.
Specifically,
	for each solver $s$, each instance $r$ and problem $p$,
	we compute the number of block evaluations $b_{s,p,r}(\tau)$ such that
	\begin{equation}
		\frac{|f_{s,b,p,r} - f_p^*|}{|f_p^*|} \le \tau,
		\label{eq:condition-perf-profile}
	\end{equation}
	where $f_p^*$ is the best known objective value for problem $p$
	and $f_{s,b,p,r}$ is the value obtained with the solver $s$
	after $b$ block evaluations.
Let $b^{\min}_{p,r}(\tau)$ be the smallest budget
	for solving the instance $r$ of problem $p$ with deviation $\tau$, i.e.,
	$$ b^{\min}_{p,r}(\tau) = \underset{s}{\min} \;  b_{s,p,r}(\tau). $$
Then, we can plot the proportion of runs of solver $s$
that satisfy Eq.~\eqref{eq:condition-perf-profile}
	at a multiple $\alpha$ of the smallest budget, i.e.,
	$\alpha \, b^{\min}_{p,r}(\tau)$ block evaluations.
 Figure~\ref{fig:perf-profile} depicts the performance profiles of the five considered problems
	for $q = 64$ over the 50 runs (instances) for $\tau$ values that range from 10$^{-1}$ to 10$^{-4}$.

\renewcommand{\figWidth}{0.45\linewidth}
\newcommand{\dir}{/figs/PerfProfiles_64}
\begin{figure}[ht!]
	\center
	\begin{tabular}{c c}
		Performance profiles for $\tau=10^{-1}$				&	 Performance profiles for $\tau=10^{-2}$ 				\\
		\includegraphics[width=\figWidth]{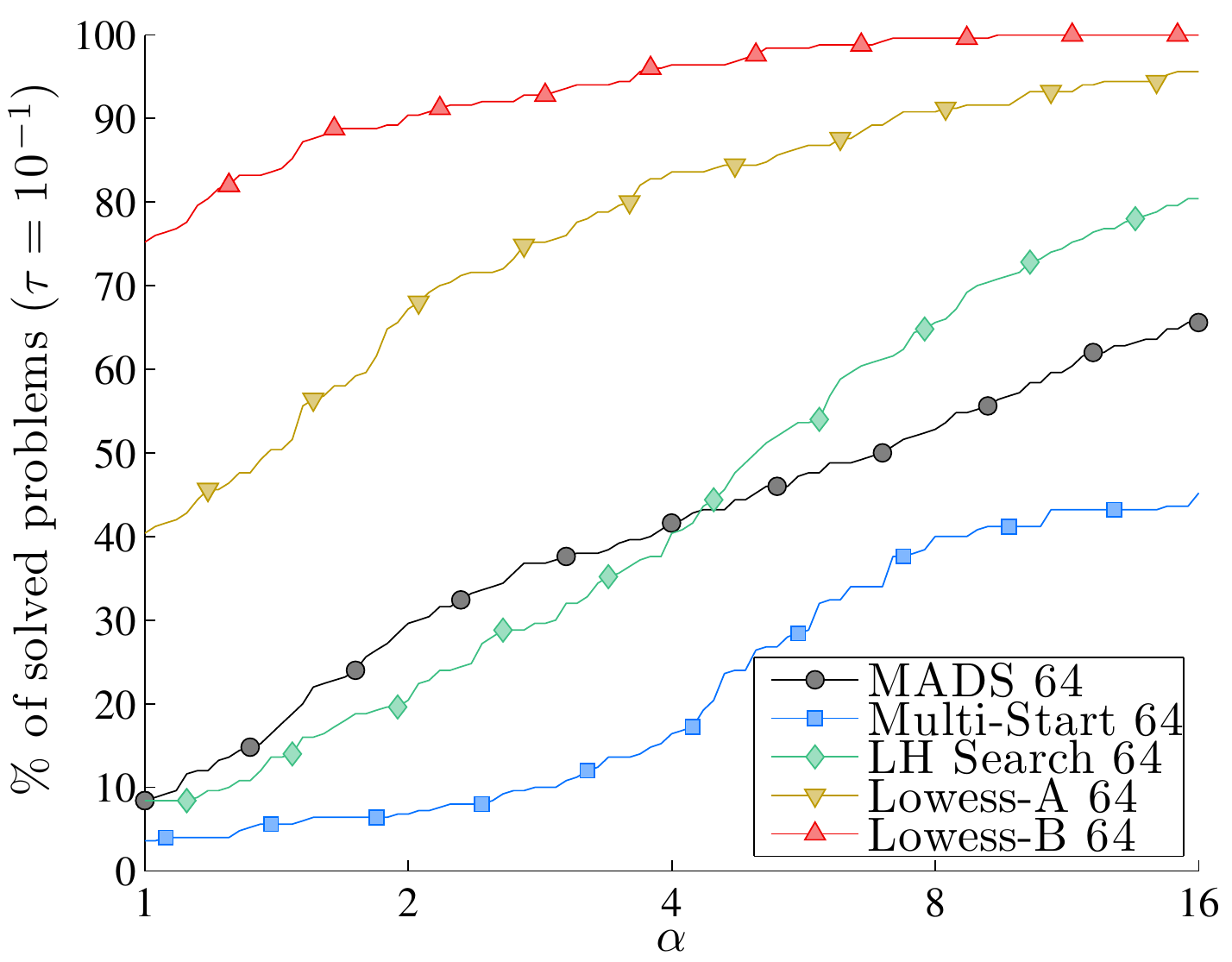}	&	\includegraphics[width=\figWidth]{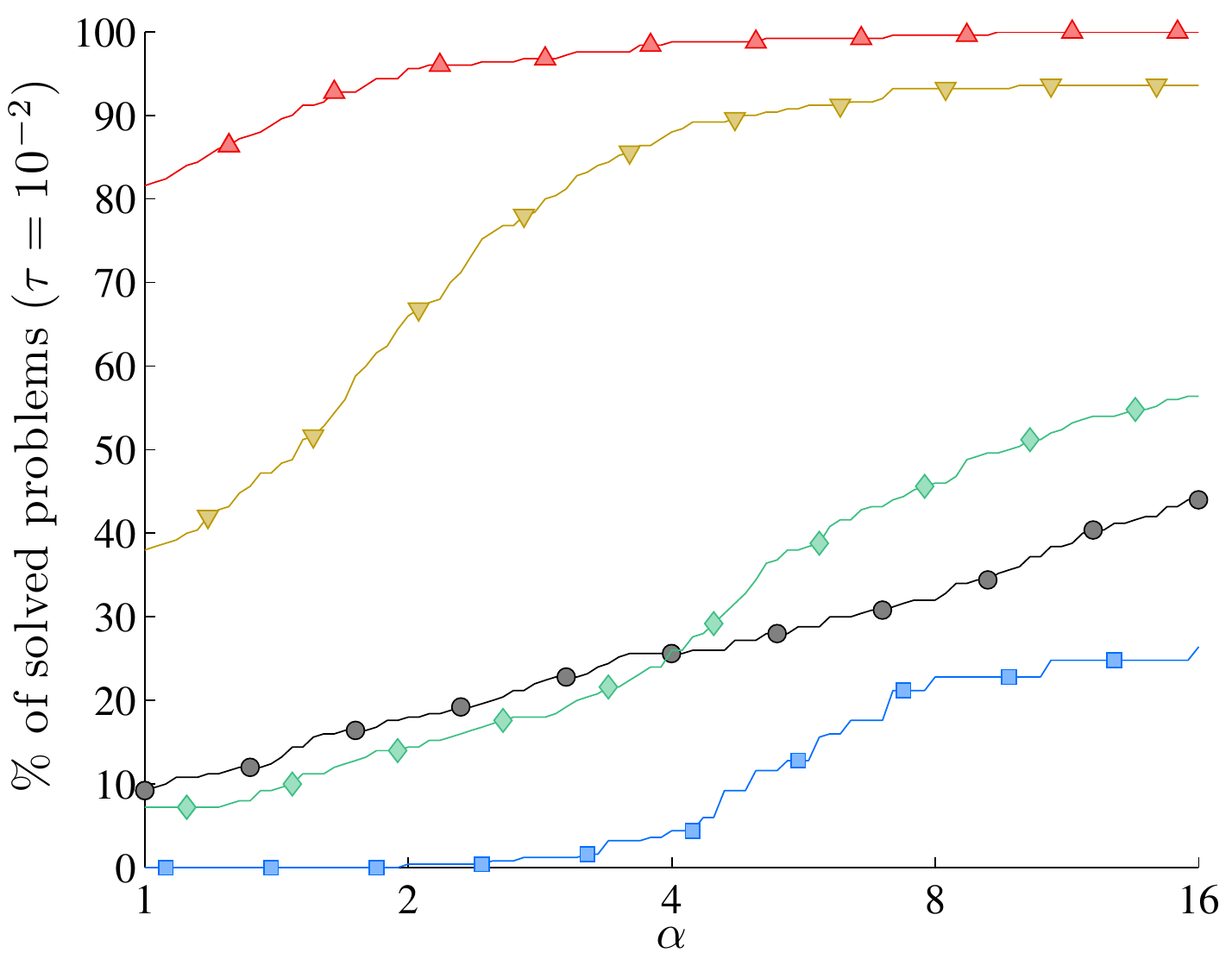}	\\
												&											\\
		Performance profiles for $\tau=10^{-3}$ 				&	 Performance profiles for $\tau=10^{-4}$				\\
		\includegraphics[width=\figWidth]{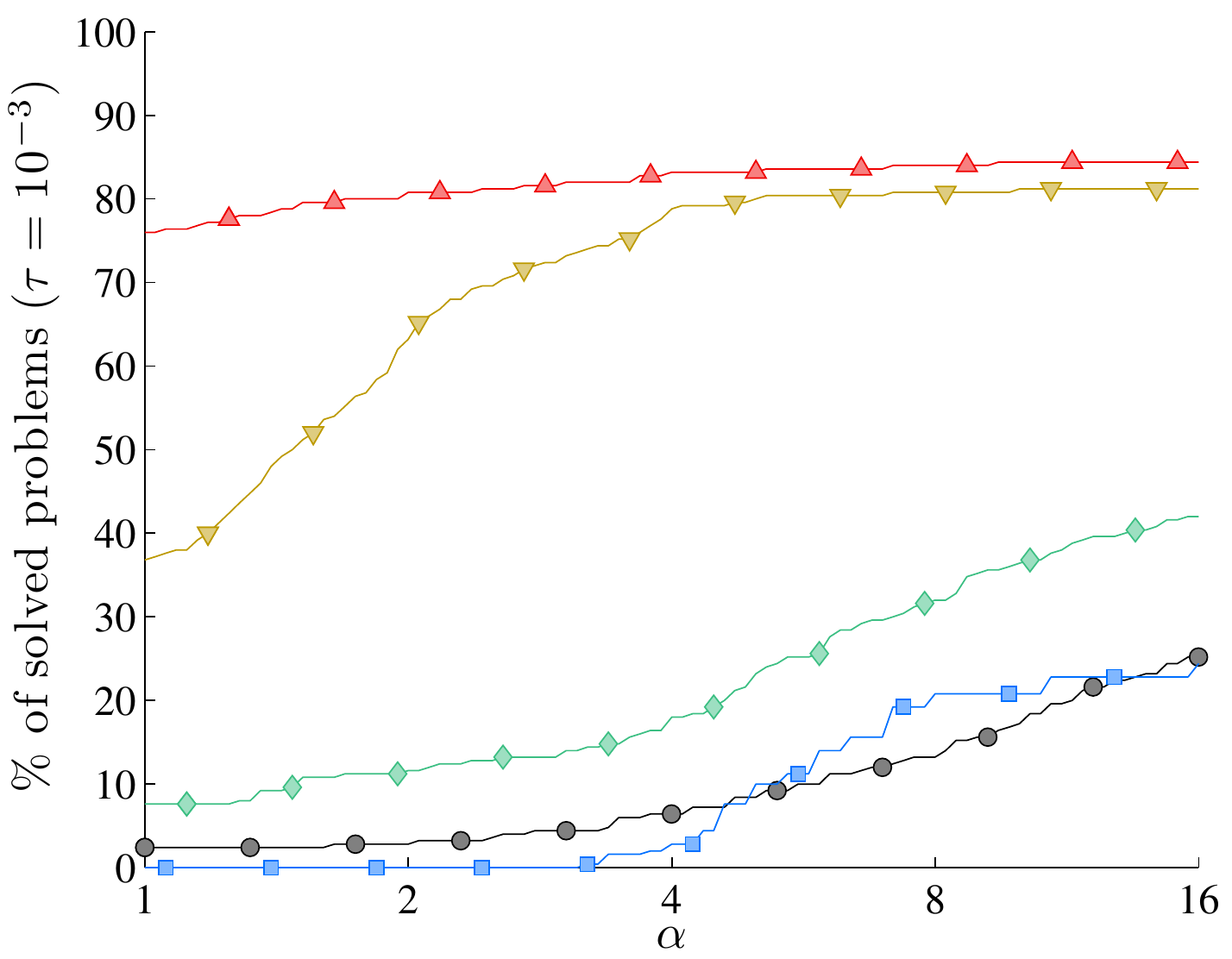} 	&	\includegraphics[width=\figWidth]{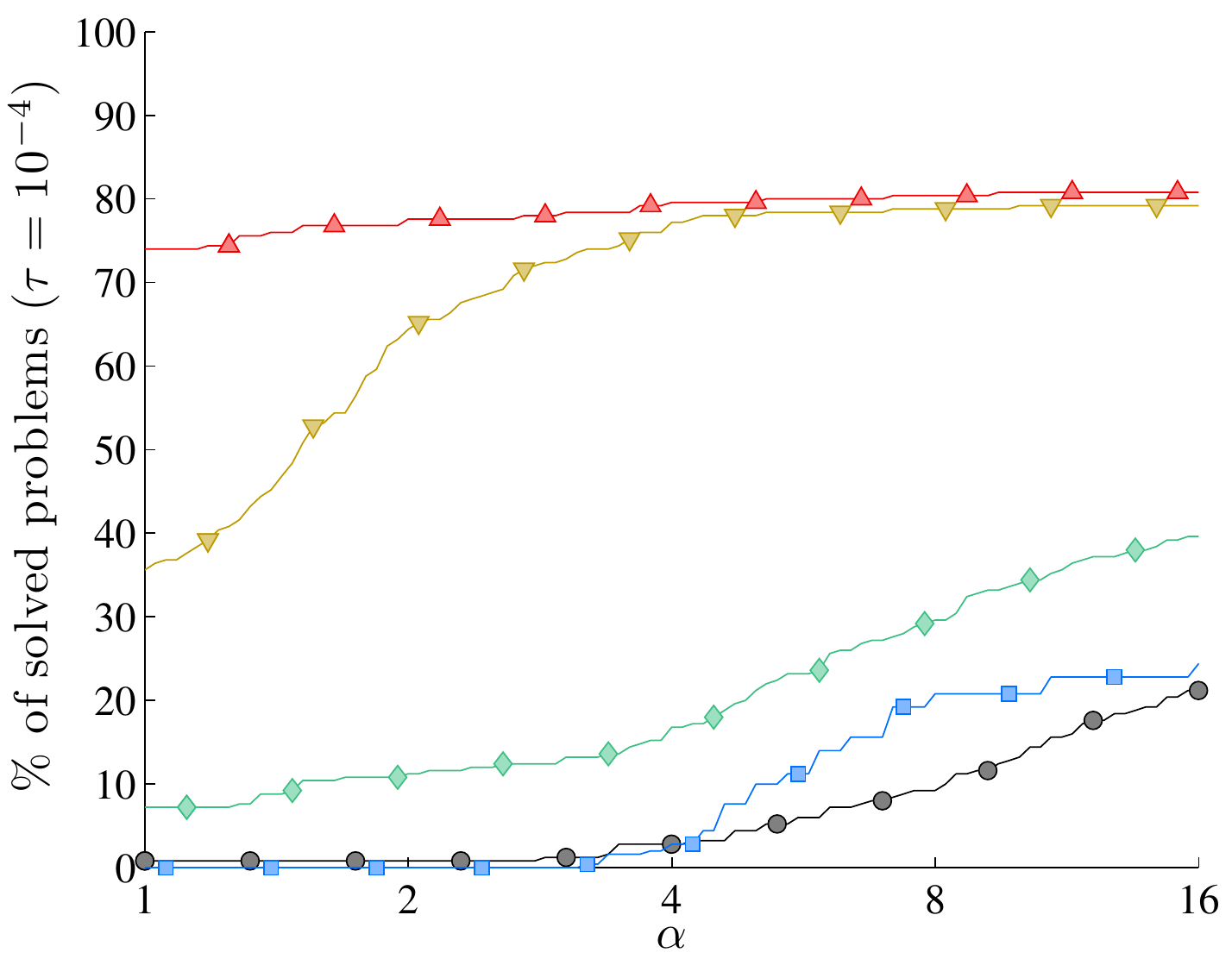}	\\
	\end{tabular}
	\caption{Performance profiles for $q=64$ over the 50 runs of all five problems}
	\label{fig:perf-profile}
\end{figure}

Note that higher curves imply better solver performance.
Moreover,
	a solver that performs well for small values of $\alpha$
	is a solver that can solve,
	for the considered $\tau$,
	a large number of problems with a small evaluation budget.
 Figure~\ref{fig:perf-profile} confirms our previous observations:
 \quote{Lowess-B} outperforms all solvers,
	followed by \quote{Lowess-A} and \quote{LH Search}
	in second and third position.
	\quote{MADS} and \quote{Multi-Start} are in the last position.
For large tolerances,
	\quote{MADS} does better than \quote{Multi-Start},
	and better than \quote{LH~Search} for small values of $\alpha$ ($\le 4$).
For small precisions,
	\quote{MADS} is outperformed by the other solvers.

The most interesting observation is the significant gap
	between the performance curves of LOWESS solvers
	and the ones from the three other solvers.
The gap increases as $\tau$ decreases.
From moderate to lower $\tau$ values ($\le 10^{-2}$),
	\quote{Lowess-B} systematically solves at least twice more problems
	than \quote{LH~Search}, \quote{Multi-Start} and \quote{MADS}.
It is unusual to observe such clear differences on performance profiles.
\quote{Lowess-A} performs almost as well as \quote{Lowess-B}, 
	particularly when $\tau$ is small
	(around $10^{-4}$),
	but needs at least  four times more block evaluations
	to achieve this.

\subsubsection*{Scalability analysis}

\begin{figure}[h]
	\center
	\includegraphics[width=\figWidth]{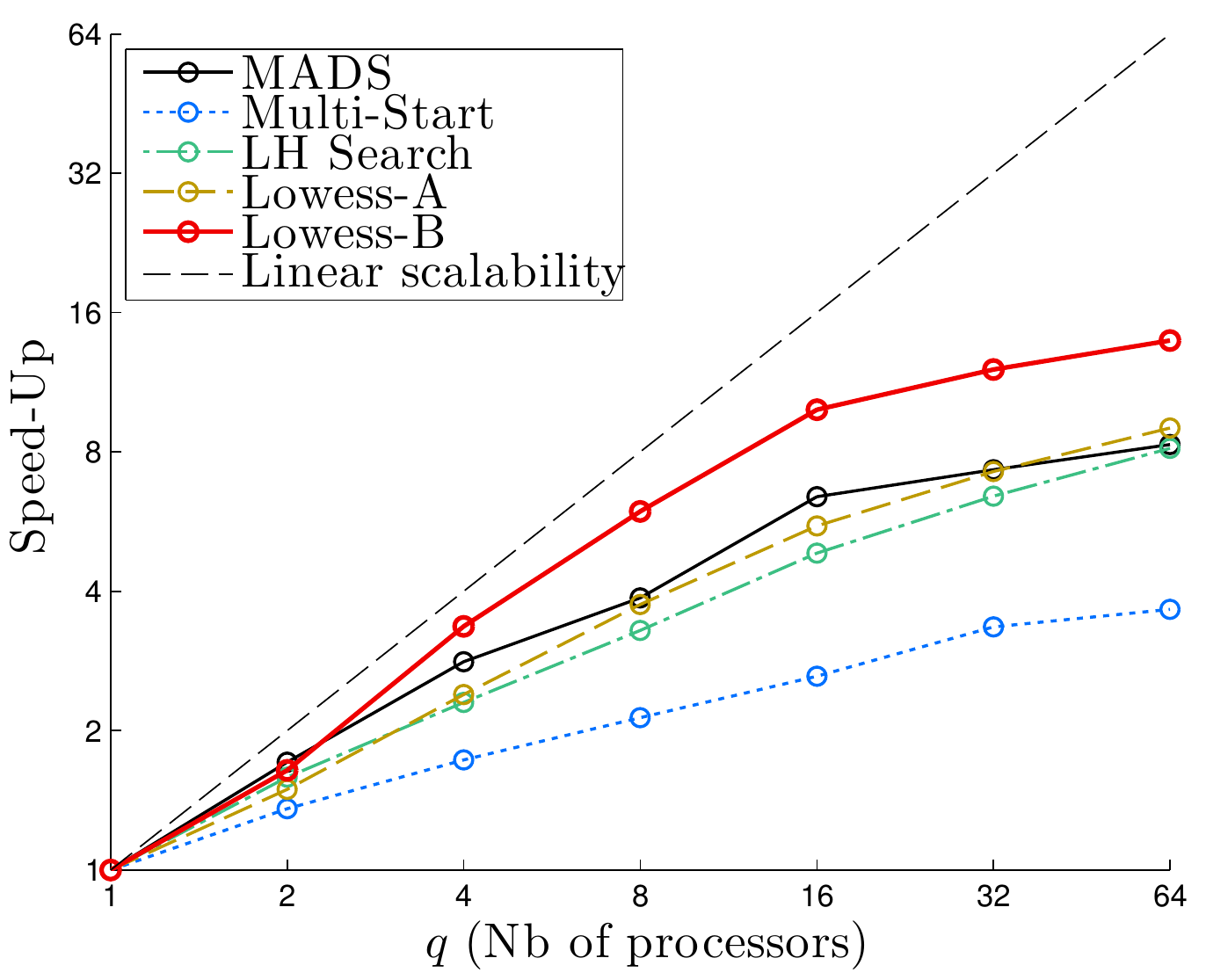}
	\includegraphics[width=\figWidth]{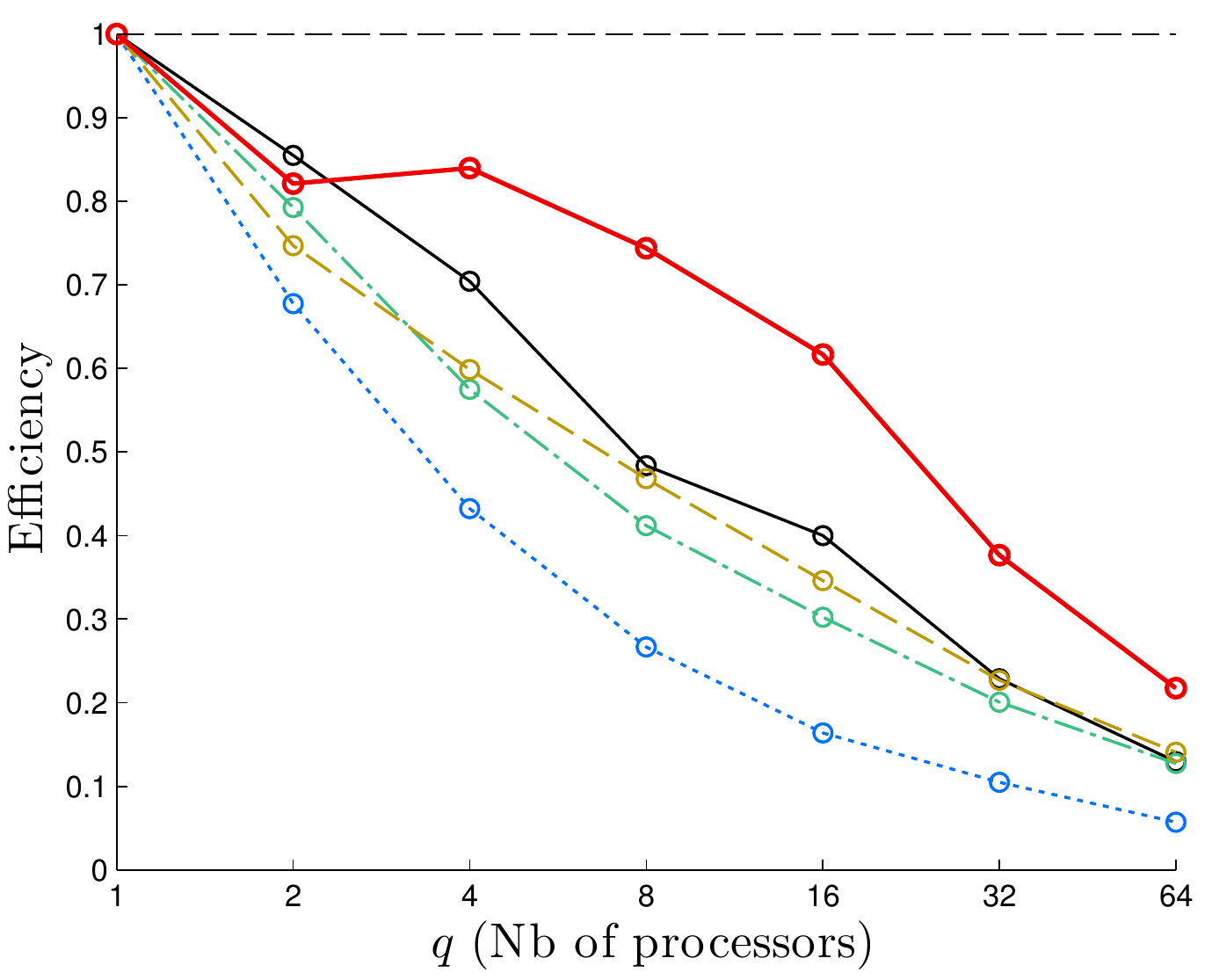}
	\caption{Speed-up and efficiency}
	\label{fig:scalability}
\end{figure}

We wish to establish the reduction of wall-clock time when using additional resources for each solver.
To that end,
	we follow the methodology proposed in~\cite{Jogalekar2000}.
We define $f_{s,b,p,r,q}$ as the value of the objective function obtained by solver $s$
	after $b$ block evaluations 
	on instance $r$ of problem $p$
	when using $q$ CPUs.
	We also define the \emph{reference} objective value as the best value achieved with only one CPU
	($q=1$), i.e.,
	\newcommand{\fref}{f^{\text{ref}}}
	\newcommand{\bref}{b^{\text{ref}}}
	\begin{equation}
		\fref_{s,p,r} = \underset{b}{\min} f_{s,b,p,r,1},
	\end{equation}
	and $\bref_{s,p,r,q}$  the number of block evaluations necessary to reach $\fref_{s,p,r}$
	when $q$ CPUs are used,
	i.e.,
	\begin{equation}
		\bref_{s,p,r,q} = \min \{ b : \fref_{s,p,r} \le f_{s,b,p,r,q} \}.
	\end{equation}
The \emph{speed-up} of solver $s$ when solving with  $q$  CPUs is defined as
	\begin{equation}
		\text{speed-up}(s,q) = \underset{p,r}{\text{geomean}} 
		\left(\;
			\frac{\bref_{s,p,r,q}}{\bref_{s,p,r,1}} 
		\;\right)
	\end{equation}
and its \emph{efficiency} as
	\begin{equation}
		\text{efficiency}(s,q) = \frac{\text{speed-up}(s,q)}{q}.
	\end{equation}
Figure~\ref{fig:scalability} depicts the speed-up and efficiency values obtained by our numerical experiments.

Perfect scalability is obtained when
	the speed-up is equal to $q$ and the efficiency is equal to 1.
The speed-up curves show that the power introduced by new CPUs decreases
	as their number increases.
This was observed in Figures~\ref{fig:curves-median-1} and~\ref{fig:curves-median-2}
	where problems exhibited saturation around $q = 8$ and 16.
\quote{Lowess-B} achieves the best speed-up,
	followed by \quote{MADS}. 
\quote{Lowess-A} and \quote{LH Search} come next, followed by the worst-performing solver, namely \quote{Multi-Start}.
	We conclude that it is better and more productive to proceed with one search performing $q$ parallel evalutions instead of conducting $q$ independent searches consisting of a single evaluation.

The efficiency curves demonstrate rapid decrease 
	except for \quote{Lowess-B}; its rate exhibits a bump on its efficiency curve at  $q = 4$.
For $q = 2$,
	only methods~\ref{algo:method3} and \ref{algo:method4} are used to generate candidates.
	For $q \geq 4$, methods~\ref{algo:method5} and \ref{algo:method6} are also used.
	The aforementioned bump highlights the important contributions of these methods to the efficiency of \quote{Lowess-B}.

\section{Conclusion}

Linear LOWESS models with optimized kernel shapes and coefficients
seem to provide high-performing surrogates of the blackboxes.
The use of diverse selection methods
	(\ref{algo:method3} to \ref{algo:method6})
	enables an efficient exploration of the design space, accelerates local convergence, and makes optimal use of additional CPU resources.
Methods~\ref{algo:method5} and \ref{algo:method6} are particularly efficient,
	outperforming the other selection methods.
This means that the way surrogates are used by method~\ref{algo:method1}
	is not effective.
Similarly, the diversification strategy of method~\ref{algo:method2}
	is not adequate to select points that lie far enough from the ones already evaluated.
	
We cannot draw a definite conclusion about which of the methods~\ref{algo:method5} or \ref{algo:method6}
	is better than the other.
We believe that the good performance of \quote{Lowess-B} is due to using method~\ref{algo:method5}. 
%

The proposed selection methods are not specific to the LOWESS model considered here; they are applicable to any surrogates. We believe that they will work well with reduced-fidelity (or variable-fidelity) physical-based models
	since high- and low-fidelity models typically have similarly structured solution domains.
The selection methods are also applicable to other algorithms using surrogates
	to identify promising points to evaluate.

\bibliographystyle{unsrt}
\bibliography{bibliography}

\begin{thebibliography}{10}

\bibitem{DBLP:series/sci/KramerCK11}
O.~Kramer, D.E. Ciaurri, and S.~Koziel.
\newblock Derivative-free optimization.
\newblock In S.~Koziel and XS. Yang, editors, {\em Computational Optimization,
  Methods and Algorithms}, volume 356 of {\em Studies in Computational
  Intelligence}, pages 61--83. Springer, 2011.

\bibitem{AuCo04a}
C.~Audet, S.~{Le~Digabel}, C.~Tribes, and V.~Rochon Montplaisir.
\newblock The {NOMAD} project.
\newblock Software available at \url{https://www.gerad.ca/nomad}.

\bibitem{Le09b}
S.~{Le~Digabel}.
\newblock Algorithm 909: {NOMAD}: Nonlinear optimization with the {MADS}
  algorithm.
\newblock {\em {ACM} Transactions on Mathematical Software}, 37(4):44:1--44:15,
  2011.

\bibitem{AuDe2006}
C.~Audet and J.E. {Dennis, Jr.}
\newblock Mesh adaptive direct search algorithms for constrained optimization.
\newblock {\em SIAM Journal on Optimization}, 17(1):188--217, 2006.

\bibitem{AlbaLuqueNesmachnow2013}
E.~Alba, G.~Luque, and S.~Nesmachnow.
\newblock Parallel metaheuristics: recent advances and new trends.
\newblock {\em International Transactions in Operational Research},
  20(1):1--48, 2013.

\bibitem{LeAbAuDe10}
S.~{Le~Digabel}, M.A. Abramson, C.~Audet, and J.E. {Dennis, Jr.}
\newblock Parallel versions of the {MADS} algorithm for black-box optimization.
\newblock In {\em Optimization days}, Montreal, May 2010. GERAD.
\newblock Slides available at
  \url{http://www.gerad.ca/Sebastien.Le.Digabel/talks/2010_JOPT_25mins.pdf}.

\bibitem{TaLeDKo2014}
B.~Talgorn, S.~{Le~Digabel}, and M.~Kokkolaras.
\newblock {Statistical Surrogate Formulations for Simulation-Based Design
  Optimization}.
\newblock {\em Journal of Mechanical Design}, 137(2):021405--1--021405--18,
  2015.

\bibitem{PoTaHaKoLe2014}
Mahdi Pourbagian, Bastien Talgorn, WagdiG. Habashi, Michael Kokkolaras, and
  S\'ebastien Le~Digabel.
\newblock Constrained problem formulations for power optimization of aircraft
  electro-thermal anti-icing systems.
\newblock {\em Optimization and Engineering}, pages 1--31, 2015.

\bibitem{Ensemble2016}
C.~Audet, M.~Kokkolaras, S.~{Le~Digabel}, and B.~Talgorn.
\newblock {Order-based error for managing ensembles of surrogates in
  derivative-free optimization}.
\newblock {\em Journal of Global Optimization}, 70(3):645--675, 2018.

\bibitem{Lowess2016}
B.~Talgorn, C.~Audet, M.~Kokkolaras, and S.~{Le~Digabel}.
\newblock {Locally weighted regression models for surrogate-assisted design
  optimization}.
\newblock {\em Optimization and Engineering}, 19(1):213--238, 2018.

\bibitem{Haftka2016}
Raphael~T. Haftka, Diane Villanueva, and Anirban Chaudhuri.
\newblock Parallel surrogate-assisted global optimization with expensive
  functions -- a survey.
\newblock {\em Structural and Multidisciplinary Optimization}, 54(1):3--13, Jul
  2016.

\bibitem{FlLe02a}
R.~Fletcher and S.~Leyffer.
\newblock Nonlinear programming without a penalty function.
\newblock {\em Mathematical Programming}, Series A, 91:239--269, 2002.

\bibitem{AuDe09a}
C.~Audet and J.E. {Dennis, Jr.}
\newblock A progressive barrier for derivative-free nonlinear programming.
\newblock {\em SIAM Journal on Optimization}, 20(1):445--472, 2009.

\bibitem{AlizadehAllenMistree2020}
R.~Alizadeh, J.K. Allen, and F.~Mistree.
\newblock Managing computational complexity using surrogate models: a critical
  review.
\newblock {\em Research in Engineering Design}, 2020.

\bibitem{Cle1979}
W.S. Cleveland.
\newblock Robust locally weighted regression and smoothing scatterplots.
\newblock {\em Journal of the American Statistical Association}, 74:829--836,
  1979.

\bibitem{Cle1981}
W.S. Cleveland.
\newblock {LOWESS}: {A} {P}rogram for {S}moothing {S}catterplots by {R}obust
  {L}ocally {W}eighted {R}egression.
\newblock {\em The American Statistician}, 35(1), 1981.

\bibitem{Cle1988a}
W.S. Cleveland and S.J. Devlin.
\newblock Locally weighted regression: {A}n approach to regression analysis by
  local fitting.
\newblock {\em Journal of the American Statistical Association}, 83:596--610,
  1988.

\bibitem{Cle1988b}
W.S. Cleveland, S.J. Devlin, and E.~Grosse.
\newblock Regression by local fitting: methods, properties, and computational
  algorithms.
\newblock {\em Journal of Econometrics}, 37(1):87 -- 114, 1988.

\bibitem{sgtelib}
B.~Talgorn.
\newblock {SGTELIB}: Surrogate model library for derivative-free optimization.
\newblock \url{https://github.com/bbopt/sgtelib}, 2019.

\bibitem{McCoBe79a}
M.D. McKay, R.J. Beckman, and W.J. Conover.
\newblock A comparison of three methods for selecting values of input variables
  in the analysis of output from a computer code.
\newblock {\em Technometrics}, 21(2):239--245, 1979.

\bibitem{SaWiNoLH2003}
T.J. Santner, B.J. Williams, and W.I. Notz.
\newblock {\em The Design and Analysis of Computer Experiments}, chapter 5.2.2,
  Designs Generated by Latin Hypercube Sampling, pages 127--132.
\newblock Springer, New York, NY, 2003.

\bibitem{MlHa97a}
N.~Mladenovi\'c and P.~Hansen.
\newblock Variable neighborhood search.
\newblock {\em Computers and Operations Research}, 24(11):1097--1100, 1997.

\bibitem{HaMl01a}
P.~Hansen and N.~Mladenovi\'c.
\newblock Variable neighborhood search: principles and applications.
\newblock {\em European Journal of Operational Research}, 130(3):449--467,
  2001.

\bibitem{Garg2014}
H.~Garg.
\newblock Solving structural engineering design optimization problems using an
  artificial bee colony algorithm.
\newblock {\em Journal of Industrial and Management Optimization},
  10(3):777--794, 2014.

\bibitem{Arora2004}
J.~Arora.
\newblock {\em Introduction to Optimum Design}.
\newblock Elsevier Science, 2004.

\bibitem{Belegundu1982}
A.D. Belegundu.
\newblock {\em A Study of Mathematical Programming Methods for Structural
  Optimization}.
\newblock University of Iowa, 1982.

\bibitem{KaKr1993}
B.~K. Kannan and S.~N. Kramer.
\newblock {Augmented Lagrange multiplier based method for mixed integer
  discrete continuous optimization and its applications to mechanical design}.
\newblock {\em Journal of Mechanical Design}, 65:103--112+, 1993.

\bibitem{rao1996}
Singiresu~S. Rao.
\newblock {\em {Engineering Optimization: Theory and Practice, 3rd Edition}}.
\newblock Wiley-Interscience, 1996.

\bibitem{MScMLG}
Mathieu~Lemyre Garneau.
\newblock {Modelling of a solar thermal power plant for benchmarking blackbox
  optimization solvers}.
\newblock Master's thesis, \'Ecole Polytechnique de Montr\'eal, 2015.

\bibitem{DoMo02}
E.D. Dolan and J.J. Mor\'e.
\newblock Benchmarking optimization software with performance profiles.
\newblock {\em Mathematical Programming}, 91(2):201--213, 2002.

\bibitem{Jogalekar2000}
Prasad Jogalekar and Murray Woodside.
\newblock Evaluating the scalability of distributed systems.
\newblock {\em IEEE Trans. Parallel Distrib. Syst.}, 11(6):589–603, June
  2000.

\bibitem{Atkeson1997}
C.G. Atkeson, A.W. Moore, and S.~Schaal.
\newblock Locally weighted learning.
\newblock {\em Artificial Intelligence Review}, pages 11--73, 1997.

\end{thebibliography}

\clearpage
\appendix

\clearpage
\section{LOWESS predictions}
\label{sec:Lowess}

As a convention,
	we denote with $\xib \in \chi \subseteq \R^n$ the point of the design space
	where we want to predict the value of the blackbox output.
Locally weighted scatterplot smoothing (LOWESS) models build a local linear regression at the point $\xib$
	where the blackbox output $[f \; c_1  \hdots c_m]$ are
	to be estimated~\cite{Atkeson1997,Cle1979,Cle1981,Cle1988a,Cle1988b,Lowess2016}.
This local regression emphasizes data points that are close to $\xib$.
The interested reader can refer to~\cite{Lowess2016} for details about the method described below.
We consider here 
	only local linear regressions; 
local quadratic regressions and Tikhonov regularization are considered in~\cite{Lowess2016}.
On the contrary, while only a Gaussian kernel was considered in~\cite{Lowess2016},
six others are added here as kernel functions.

We define 
	the output matrix $\Y \in \R^{p\times (m+1)}$, 
	the design matrix $\Z_\xib \in \R^{p\times (n+1)}$, and
	the weight matrix $\W_\xib \in \R^{p\times p}$: 
	\begin{equation}
		\Y = \left[
			\begin{array}{c c @{\,} c @{\,} c}
				f(\x_1)	& c_1(\x_1)	& \hdots	& c_m(\x_1)	\\
				\vdots 	&  \vdots 	&		& \vdots 	\\
				f(\x_p)	& c_1(\x_p)	& \hdots	& c_m(\x_p)
			\end{array}
		\right],
		\;\;
		\Z_\xib = \left[
			\begin{array}{c c}
				1 		& (\x_1-\xib)\t 	\\
				\vdots 	& \vdots 		\\
				1 		& (\x_p-\xib)\t 
			\end{array}
		\right],
		\;\;
		\W_\xib = \left[
			\begin{array}{c @{} c @{} c}
				w_1(\xib) & & \\
				& \ddots &  \\
				& & w_p(\xib)\t 
			\end{array}
		\right].
	\end{equation}
The details of the computation of $w_i(\xib)$ are described in Section~\ref{sec:weights}.
Then,
	we define $\u_\xib \in \R^{n+1}$ as the first column of $(\Z_\xib^\top \W_\xib \Z_\xib)^{-1}$,
	which means that $\u_\xib$ is the solution of the linear system $\Z_\xib^\top \W_\xib \Z_\xib \u_\xib = \e_1$.
The prediction of the blackbox outputs at $\xib$ is then
	\begin{equation}
		\hat{\y}(\xib) = \left[
			\begin{array}{c c c c}
				\fh(\xib) & \ch_1(\xib) & \hdots & \ch_m(\xib) 
			\end{array}
		\right]
		= \u_\xib\t\Z_\xib^\top \W_\xib \Y.
	\end{equation}
The cross-validation value $\ycv(\x_i)$
	(i.e., the value of the LOWESS model at $\x_i$ when the data point $\x_i$ is not used to build the model)
	are computed by setting $w_i$ to 0.
Unfortunately,
	unlike for radial basis function (RBF) models or polynomial response surfaces (PRSs),
	we do not know any computational shortcut allowing a more efficient computation of the values of $\ycv$.
However,
	each value $\ycv(\x_i)$ is computed at the same computational cost as a prediction $\yh(\x_i)$.

\subsection{Weights computation in LOWESS models}
\label{sec:weights}

\begin{figure}[ht!]
	\center
	\includegraphics[width=9cm]{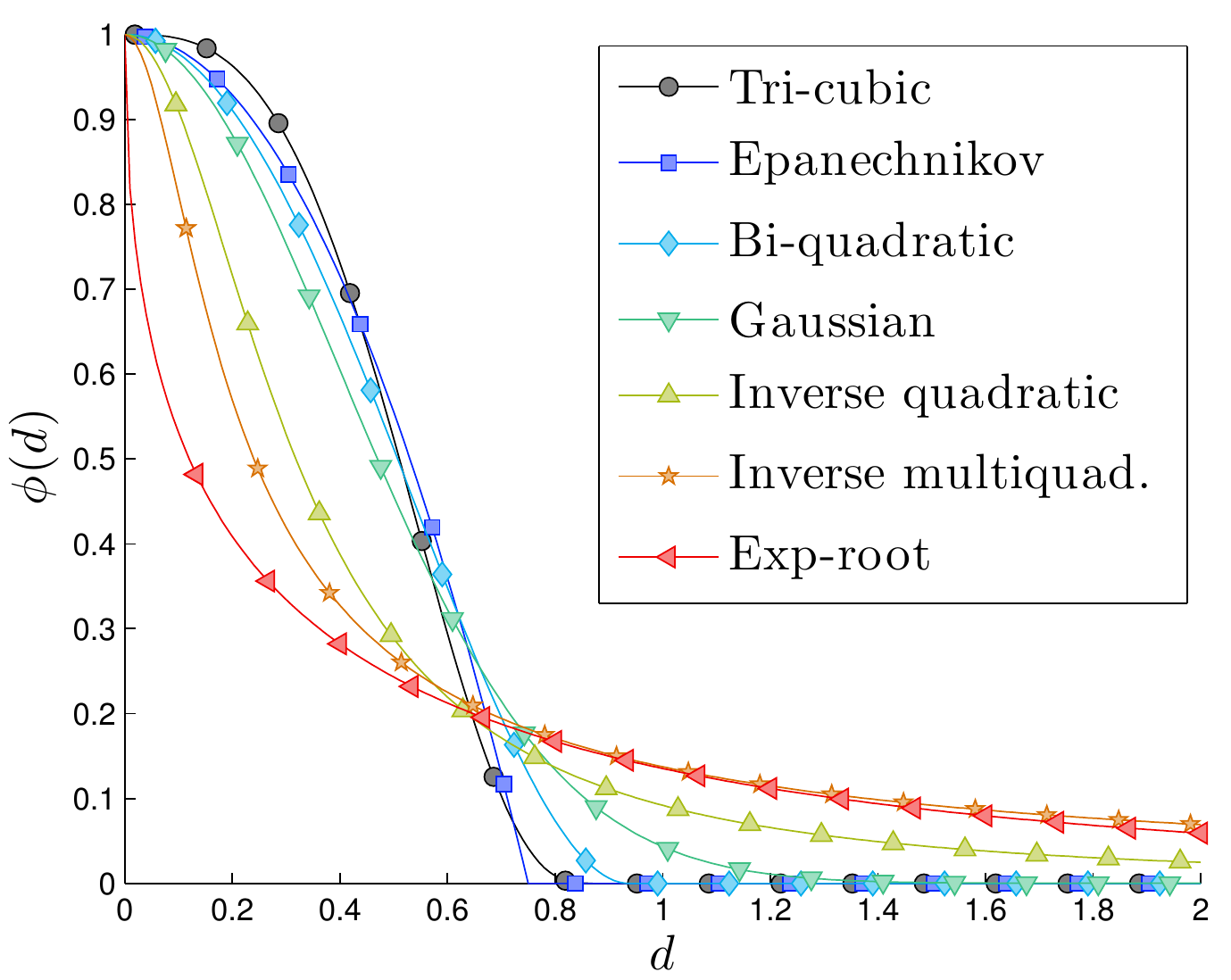}
	\caption{Representation of the 7 kernels listed in Table~\ref{tab:kernel}}
	\label{fig:kernel}
\end{figure}

The weight $w_i(\xib)$ quantifies the relative importance of the data point $\x_i$
	in the construction of the local regression at $\xib$.
Like for kernel smoothing,
	it relies on a kernel function $\phi$ and depends on the distance between $\xib$ and $\x_i$.
In our method,
	we use
	\begin{equation}
		w_i(\xib) = \phi \left (
			\lambda \frac{\|\xib-\x_i\|_2}{d_{n+1}(\xib)}
		\right),
	\end{equation}
	where $\phi(d)$ is one of the kernel functions described in Table~\ref{tab:kernel} and Figure~\ref{fig:kernel}.
All kernel functions are normalized so that $\phi(0) = 1$ and,
	if applicable,
	$\int_\R\phi = 1$.
As the integral of the inverse multi-quadratic kernel does not converge,
	the normalization constant 52.015 is introduced to minimize the $\mathcal{L}^2$ distance
	between the inverse multi-quadratic and inverse quadratic kernel.
The parameter $\lambda>0$ controls the general shape of the model,
	and $d_{n+1}(\xib)$ is a local scaling coefficient
	that estimates the distance of the $n+1^{th}$ closest data point to $\xib$.
The kernel function $\phi$ and the shape parameter $\lambda$ are chosen
	to minimize the aggregate order error with cross-validation (AOECV) described in Section \ref{sec:aoecv}.
The fact that some of the available kernel function have a compact domain
	gives to LOWESS models the ability to ignore outliers or aberrant data points.
As an example,
	if the blackbox fails to compute correctly the objective function for a given data point,
	the value returned by the blackbox might be an arbitrarily high value
	(e.g., $1.8~10^{308}$ for a C++ code returning the standard max double).
With non-compact kernel function,
	this would perturb the LOWESS model on the entire design space.
However,
	if there is no such aberrant data points,
	non-compact kernel functions tend to yield better results.

\begin{table}[ht!]
	\def\arraystretch{1.4}
	\center 
	\caption{Possible values for the kernel function $\phi$}
	\begin{tabular}{c l l c}
		\hline
		\# 	& Kernel name & $\phi:\R\rightarrow\R^+$ & Compact domain   \\
		\hline\hline
		1 	& Tri-cubic			& $\phi(d) = (1-|\frac{162}{140}d|^3)^3\ind_{|d|\le\frac{140}{162}}$	& Yes\\
		2 	& Epanechnikov		& $\phi(d) = (1-\frac{16}{9}d^2)\ind_{|d|\le\frac{3}{4}}$			& Yes\\
		3 	& Bi-quadratic		& $\phi(d) = (1-|\frac{16}{15}d|^2)^2\ind_{|d|\le\frac{15}{16}}$	& Yes\\
		4 	& Gaussian			& $\phi(d) = \exp(-\pi d^2)$							& No\\ 
		5 	& Inverse quadratic		&  $\phi(d) = \frac{1}{1+\pi^2d^2}$						& No\\
		6 	& Inverse multi-quadratic	& $\phi(d) = \frac{1}{\sqrt{1+52.015d^2}}$					& No\\
		7 	& Exp-root			& $\phi(d) = \exp(-2\sqrt{|d|})$							& No\\ 
		\hline
	\end{tabular}
		\label{tab:kernel}
	\def\arraystretch{1}
\end{table}

To obtain a model $\yh$ that is differentiable everywhere,
	\cite{Lowess2016} defines $d_{n+1}(\xib)$ such that
	the expected number of training points in a ball of center $\xib$ and radius $d_{n+1}(\xib)$ is $n+1$:
	\begin{equation}
		\E
		\left[
			\card\Big\{
				\x_i : \x_i\in\X, \|\xib-\x_i\|_2 \le d_{n+1}(\xib)
			\Big\}
		\right] = n+1.
		\label{eq:expected-card}
	\end{equation}

Moreover,
	\cite{Lowess2016} observes that the values $\big\{\|\xib-\x_i\|_2^2\big\}_{i=1,\hdots,p}$
	can be fitted well by a Gamma distribution
	and therefore defines the local scaling parameter as
	\begin{equation}
		d_{n+1}(\xib) = \sqrt{
			g^{(-1)}\left(\frac{\mu_\xib^2}{\sigma_\xib^2},\frac{\sigma_\xib^2}{\mu_\xib} ; \frac{n+1}{p} \right)
		},
	\end{equation}
	where $\mu_\xib$ (resp. $\sigma^2_\xib$) denotes the mean (resp. variance) of $\|\xib-\x_i\|_2^2$ over $\X$
	and $g^{(-1)}(k,\theta;.)$ is the inverse function of the cumulative density function of a Gamma distribution
	with shape parameter $k$ and scale parameter $\theta$.

\subsection{Aggregate Order Error with Cross-Validation}
\label{sec:aoecv}

The AOECV is an error metric that aims at quantifying the quality of a \emph{multi-output} surrogate model.
Specifically,
	it aims at quantifying the discrepancy between problems~\eqref{eq:MainProblem} and~\eqref{eq:SurrogateProblem}
	for a given surrogate model.
We first define the aggregate constraint violation function~\cite{FlLe02a}
	$h(\x) = \sum_{j=1}^m \max\{0,c_j(\x)\}^2$.
Note that other definitions of $h$ are possible
	(notably: number of violated constraints, most violated constraint, etc.)
	but as the previous definition of $h$ is used in the main MADS instance to solve~\eqref{eq:MainProblem},
	we need to use the same aggregate constraint in our definition of the AOECV.

We then define the order operators
	\begin{align}
		\x \prec \x' \Leftrightarrow &	\left\{
			\begin{array}{l}
				h(\x) < h(\x')	\\
				\text{or}		\\
				h(\x)=h(\x') \text{ and } f(\x)<f(\x'),
			\end{array}
		\right.\\
		\x \preceq \x' \Leftrightarrow & \;\; {\tt not}(\x' \prec \x)
	\end{align}
	which are transitive.
In particular,
	the incumbent solution $\x^t$ of the original problem~\eqref{eq:MainProblem}
	is such that $\x^t\preceq \x,\; \forall \x \in \X$.
Similarly,
	a global minimizer $\x^*$ is such that $\x^*\preceq \x,\; \forall \x \in \chi$.
By the same principle,
	we define the operator $\preccv$ by using the cross-validation values
	$\fcv$ and $\hcv = \sum_{j=1}^m \max\{0,\ccv_j(\x)\}^2$ instead of $f$ and $h$.
We then define the aggregated order error with cross-validation (AOECV) metric:
	\begin{equation}
		\metric_{AOECV} = 
		\frac{1}{p^2} \displaystyle \sum_{i=1}^p\sum_{j=1}^p \xor \Big(
			\x_i \prec \x_j , \x_i \preccv \x_j
		\Big).
		\label{eq:AOECV}
	\end{equation}
	where \xor is the exclusive {\tt or} operator
	(i.e., $\xor(A,B)=1$ if the booleans $A$ and $B$ differ and 0 otherwise).
The metric allows to quantify how often the model is able to correctly decide
	which of two points is better.

The shape parameter $\lambda$ and the kernel function $\phi$ are then chosen to minimize 
	$\metric_{AOECV}(\lambda,\phi).$
If two couples $(\lambda,\phi)$ lead to the same metric value
	(because of the piecewise-constant nature of the metric),
	the couple with the smallest value of $\lambda$ (i.e.,  the smoother model) is preferred.

\clearpage
\section{Detailed description of the test problems}
\label{annexe:detailedProb}

The five engineering design application problems considered
	are listed in Table~\ref{tab:5pb}.
Problem size is reflected by $n$ and $m$,
	where $n$ denotes the number of design variables
	and $m$ the number of general nonlinear inequality constraints.
Table~\ref{tab:5pb} also indicates
	whether any variables are integer or unbounded and reports the best known value of the objective function.

\begin{table}[ht!]
	\center
	\caption{Summary of the five engineering design optimization problems}
	\begin{tabular}{ l c c c c c c }
		\hline
		Problem	&  $n$		&  $m$	&   Integer	&   Infinite	&   Best objective		\\
		name		&		&		&   variables	&   bounds	&  function value	\\ 
		\hline\hline
		TCSD		&  3		&  4		&   No		&   No		&   0.0126652	\\
		Vessel		&  4		&  4		&   No		&   No		&   5,885.332	\\
		Welded	&  4		&  6		&   No		&   No		&   2.38096		\\
		Solar~1	&  9		&  5		&   Yes	&   Yes	&  --900,417		\\
		Solar~7	&  7		&  6		&   Yes	&   Yes	&  --4,976.17	\\
		\hline
	\end{tabular}
	\label{tab:5pb}
\end{table}

%
%
%
%
%



The Tension/Compression Spring Design (TCSD) problem
	consists of minimizing the weight of a spring under mechanical constraints ~\cite{Garg2014, Arora2004, Belegundu1982}.
The design variables define the geometry of the spring.
The constraints concern shear stress, surge frequency and minimum deflection.
The best known solution,
	denoted $\xstar$,
	and the bounds on the variables,
	denoted by $\xlb$ and $\xub$, 
	are given  in  Table~\ref{tab:var-tcsd}.

\begin{table}[ht!]
 	\center
	\caption{Variables of the TCSD problem}
	\begin{tabular}{ l c c c }
		\hline
		Variable description		& $\xlb$	& $\xub$	& $\xstar$ 			\\
		\hline\hline
		Mean coil diameter		& 0.05	& 2		& 0.051686696913218 	\\
		Wire diameter		& 0.25	& 1.3 		& 0.356660815351066 	\\
		Number of active coil	& 2		& 15 		& 11.292312882259289 	\\
		\hline
	\end{tabular}
	\label{tab:var-tcsd}
\end{table}


The Vessel problem
	considers the optimal design of a compressed air storage tank~\cite{Garg2014, KaKr1993}.
The  design variables define the geometry of the tank
The  constraints are related to the volume,
	pressure,
	and solidity of the tank.
The objective is to minimize the total cost of the tank,
	including material and labour.
Table~\ref{tab:var-vessel} lists the variable bounds and the best known solution.

\begin{table}[ht!]
 	\center
	\caption{Variables of the Vessel problem}
	\begin{tabular}{ l c c c }
		\hline
		Variable description 			& $\xlb$	 & $\xub$ 	& $\xstar$ 			\\
		\hline\hline
		Thickness of the vessel			& 0.0625	& 6.1875	&    0.778168641330718	\\
		Thickness of the head			& 0.0625	& 6.1875	&    0.384649162605973	\\
		Inner radius					& 10		& 200		&  40.319618721803231	\\
		Length of the vessel without heads	& 10		& 200		& 199.999999998822659	\\
		\hline
	\end{tabular}
	\label{tab:var-vessel}
\end{table}



The Welded
	(or welded beam design)
	problem (Version I)	consists of minimizing the construction cost of a beam
	under shear stress,
	bending stress,
	load
	and deflection constraints~\cite{Garg2014, rao1996}.
The  design variables define the geometry of the beam
	and the characteristics of the welded joint.
Table~\ref{tab:var-welded} lists the variable bounds and the best known solution.

\begin{table}[ht!]
	\center
	\caption{Variables of the Welded problem}
	\begin{tabular}{ l c c c }
		\hline
		Variable description		& $\xlb$	& $\xub$ 	& $\xstar$ 			\\
		\hline\hline
		Thickness of the weld	&  0.1 	& 2 		& 0.244368407428265	\\
		Length of the welded joint	&  0.1 	& 10 		& 6.217496713101864	\\
		Width of the beam		&  0.1 	& 10 		& 8.291517255567012	\\
		Thickness of the beam	&  0.1 	& 2 		& 0.244368666449562	\\
		\hline
	\end{tabular}
	\label{tab:var-welded}
\end{table}



The Solar1 and Solar7 problems consider the optimization of a solar farm,
	including the heliostat field
	and/or the receiver~\cite{MScMLG}.
The Solar1 optimization problem aims at maximizing the energy received over a period of 24 hours
	under several constraints of budget
	and heliostat field area.
This problem has one integer variable that has no upper bound.
Table~\ref{tab:var-solar1} lists the variable bounds and the best known solution.

%



\begin{table}[ht!]
	\center
\caption{Variables of the Solar1 problem}
	\begin{tabular}{ l c c c }
		\hline
		Variable description 			& $\xlb$ 	& $\xub$ 	& $\xstar$ 			\\
		\hline\hline
		Heliostat height 				& 1 		& 40 		& 6.165258994385601 	\\
		Heliostat width 				& 1 		& 40 		& 10.571794049143792 	\\
		Tower height 				& 20 		& 250 		& 91.948461670428486 	\\
		Receiver aperture height			& 1 		& 30 		& 6.056202026704944 	\\
		Receiver aperture width 			& 1 		& 30 		& 11.674984434929991	\\
		Max number of heliostats (Integer)	& 1 		& $+\infty$	& 1507 			\\
		Field maximum angular span 		& 1 		& 89 		& 51.762281627953051 	\\
		Minimum distance to tower 		& 0.5		& 20 		& 1.347318830713629 	\\
		Maximum distance to tower 		& 1 		& 20 		& 14.876940809562798	\\
		\hline
	\end{tabular}
		\label{tab:var-solar1}
\end{table}


The Solar7 problem aims at maximizing the efficiency of the receiver over a period of 24 hours,
	for a given heliostats field,
	under 6 binary constraints~\cite{MScMLG}.
This problem has one integer variable that has no upper bound.
The objective function is the energy transferred to the molten salt.
Table~\ref{tab:var-solar7} lists the variable bounds and the best known solution.

\begin{table}[ht!]
	\center
	\caption{Variables of the Solar7 problem}
	\begin{tabular}{ l c c c }
		\hline
		Variable description			 & $\xlb$ 	& $\xub$ 	& $\xstar$			\\
		\hline\hline
		Aperture height 			& 1 		& 30 		& 11.543687848308958	\\
		Aperture width 			& 1 		& 30 		& 15.244236061098078	\\
		Outlet temperature	 		& 793 		& 995 		& 803.000346734710888	\\
		Number of tubes (Integer)		& 1 		& $+\infty$	& 1292			\\
		Insulation thickness 		& 0.01 	& 5 		& 3.399190219909724	\\
		Tubes inside diameter 		& 0.005	& 0.1 		& 0.010657067457678	\\
		Tubes outside diameter 		& 0.0055 	& 0.1 		& 0.011167646941518	\\
		\hline
	\end{tabular}
	
	\label{tab:var-solar7}
\end{table}

%




\end{document}